\documentclass{gtpart}     % Basic GT/GTM/AGT style
\usepackage{graphicx}
\usepackage[all]{xy}
\usepackage{graphicx}
\usepackage[usenames,dvipsnames]{color}
\usepackage[normalem]{ulem}

\title{Approximating $C^{1,0}$-foliations}

\author{William H. Kazez}
\givenname{William}
\surname{Kazez}
\address{Department of Mathematics, University of Georgia, Athens, GA 30602}
\email{will@math.uga.edu}
\author{Rachel Roberts}
\givenname{Rachel}
\surname{Roberts}
\address{Department of Mathematics, Washington University, St. Louis, MO 63130}
\email{roberts@math.wustl.edu}

\arxivreference{math.GT/1404.5919v2}
\arxivpassword{8qa7w}

\volumenumber{}
\issuenumber{}
\publicationyear{}
\papernumber{}
\startpage{}
\endpage{}
\doi{}
\MR{}
\Zbl{}
\received{}
\revised{}
\accepted{}
\published{}
\publishedonline{}
\proposed{}
\seconded{}
\corresponding{}
\editor{}
\version{}

\newtheorem{thm}{Theorem}[section]
\newtheorem{notation}[thm]{Notation}
\newtheorem{lemma}[thm]{Lemma}
\newtheorem{cor}[thm]{Corollary}
\newtheorem{prop}[thm]{Proposition}
\newtheorem{con}[thm]{Conjecture}

\theoremstyle{definition}
\newtheorem{definition}[thm]{Definition}
\newtheorem{example}[thm]{Example}
\newtheorem{question}[thm]{Question}

\begin{document}

\begin{abstract}    % type your abstract below
We extend the Eliashberg-Thurston theorem on approximations  of taut oriented $C^2$-foliations of 3-manifolds by both positive and negative contact structures to a large class of taut oriented $C^{1,0}$-foliations,   where by $C^{1,0}$ foliation, we mean a foliation with continuous tangent plane field. These $C^{1,0}$-foliations can therefore be approximated by weakly symplectically fillable, universally tight, contact structures. This allows applications of $C^2$-foliation theory to contact topology and Floer theory to be generalized and extended to constructions of $C^{1,0}$-foliations.
\end{abstract}

\maketitle

%%%%%%%%%%%%%%%%%%%%   Start of main body of article

\section{Introduction}

 In \cite{ET}, Eliashberg and Thurston introduce the notion of {\it confoliation} and prove that when $k\ge 2$, a transversely oriented $C^k$~foliation $\mathcal F$ of a closed, oriented 3-manifold not equal to $S^1\times S^2$ can be $C^0$~approximated by a pair of $C^{\infty}$ contact structures, one positive and one negative.  They also prove that when $\mathcal F$ is also taut, any contact structure sufficiently close to the tangent plane field of $\mathcal F$ is weakly symplectically fillable and universally tight.
 For the most part, they restrict attention to confoliations which are at least $C^1$. For their main theorem, they restrict attention to confoliations which are at least $C^2$.  Their assumption that 2-plane fields and 1-forms are $C^1$ is necessary for it to be possible to take derivatives. Their assumption that 2-plane fields and 1-forms are $C^2$ is necessary for Sacksteder's Theorem (\cite{S}, see also Theorem 1.2.5 of \cite{ET}) to apply. A weakening of this $C^2$ assumption in the neighborhood of compact leaves is used by Kronheimer, Mrowka, Ozsv\'ath, and Szab\'o in \cite{KMOS} to show that the methods of \cite{ET} apply to those foliations constructed by Gabai in \cite{G1,G2,G3} which are $C^{\infty}$ except  along torus leaves.  In this paper we show that many of the techniques of Eliashberg and Thurston extend to transversely oriented, taut, $C^{1,0}$-foliations satisfying a natural transitivity condition.

Eliashberg and Thurston's proof that sufficiently smooth transversely oriented taut foliations of 3-manifolds can be perturbed to weakly symplectically fillable contact structures gives a direct connection between foliation theory and symplectic topology via contact topology.  This connection has been most spectacularly exploited by Kronheimer and Mrowka \cite{KrMr} in their proof of the Property~P Conjecture.

Since the contact structures produced by Eliashberg-Thurston are  weakly symplectically fillable and universally tight, their theorem is important in contact topology. To apply their theorem, of course one must start with a taut oriented foliation.  The issue that we seek to address is that most new constructions of foliations produce foliations that are not smooth enough to apply the Eliashberg-Thurston theorem.

Constructions of foliations which fail even to be $C^1$ can be found in \cite{DL, DR, G1, g1, g2, g3, G2, G3, Ga, KRo, KR2, Li, Li2, LR, R, R1, R2}.  In each of these papers, foliations are constructed using branched surfaces.  Foliations constructed using branched surfaces will have smooth leaves but will often vary only continuously in a transverse direction because of the role and nonsmooth nature of leaf blow-up, Operation 2.1.1 of  \cite{Ga}. (To obtain the smooth foliations  found in  \cite{G1,g1,g2,g3,G2,G3}, Gabai   avoids appealing to Operation 2.1.1 and instead takes  advantage of the fact that the branched surfaces involved in his construction are finite depth.)

The main results of this paper are, roughly, that $C^{1,0}$-foliations satisfying a simple transitivity condition can be approximated by weakly symplectically fillable contact structures. Such foliations include most of  those found in \cite{DL, DR, g1, g2, g3, G3, Ga, KRo, KR2, LR, R, R1, R2}  and many of those found in \cite{G1, G2,Li, Li2}.  In proving this, we give new methods for both initiating and propagating contact structures when beginning with a foliation.  Before giving a more precise indication of our results, we give a brief description of some of the main ideas in Eliashberg and Thurston's proof.

Eliashberg and Thurston interpolate between the notions of foliation and contact structure by introducing confoliations.  This structure restricts to a contact structure on an open set, but is equal to the original foliation elsewhere.  Their first step is to create a contact zone in the manifold, that is, to approximate the foliation so that it is a contact structure on a non-empty set.  One place they accomplish this is in a neighborhood of a curve with attracting holonomy.  We obtain a related result, Theorem~\ref{attracting}, for the larger class of curves for which the holonomy has  a contracting interval. We also introduce the notion of   $L$-bracketed foliation in Definition~\ref{transtranv} and show how a contact zone can be naturally introduced about a regular neighborhood of the corresponding link $L$.  This has applications to spines of open book decompositions, and more generally to manifolds obtained by surgery or filling.

Their next step is to propagate the contact structure throughout the manifold using smooth foliation charts.  See Lemma~\ref{ETlemma}, and Corollary~\ref{approximate}.  We use a similar strategy, but we rely on the local existence of smooth approximations to a given $C^{1,0}$-foliation.

To prove that the contact structure they produce is weakly symplectically fillable, they use the crucial structure of a smooth volume preserving flow $\Phi$ transverse to the taut foliation.  Such a smooth flow is the starting point for our work, for they exist even for $C^{1,0}$-foliations.

There are several notions that describe the relationship between a foliation $\mathcal F$ and ``nearby'' contact structures $\xi$.  If there is a continuous family of contact structures $\xi_t$ such that $\xi_t \to \mathcal F$, we say $\mathcal F$ can be {\sl perturbed}  or {\sl deformed} to a contact structure. A weaker notion is the existence of a sequence, $\xi_n \to \mathcal F$, in which case $\mathcal F$ can be {\sl approximated} by contact structures. Both notions of convergence can be refined by defining convergence using a $C^k$ norm on tangent planes for values of $k$ ranging from 0 to $\infty$.  We use the flow $\Phi$ transverse to $\mathcal F$ to define a weaker, topological notion of approximation.

\begin{definition}\label{Phiapprox} Given an oriented foliation $\mathcal F$ and a positively transverse flow $\Phi$, we say an oriented 2-plane field $\chi$, typically a contact structure or a confoliation, is a {\sl $\Phi$-approximation} of $\mathcal F$ if $\chi$ is positively transverse to $\Phi$.  We say that two oriented 2-plane fields are {\sl $\Phi$-close} if both are positively transverse to $\Phi$.
\end{definition}

One application of our work is to complete the proof of a theorem of \cite{HKM2}.  Honda, Kazez, and Mati\'c show that sufficiently large fractional Dehn twisting for an open book decomposition with connected binding implies that the canonically associated contact structure is weakly symplectically fillable.  Their proof requires the existence of contact structures approximating $C^{1,0}$-foliations constructed by Roberts in \cite{R2}, and thus needs  a stronger version of the Eliashberg-Thurston theorem.

Another application of our work is to prove that 3-manifolds containing taut, oriented $C^{1,0}$-foliations satisfying our transitivity condition cannot be L-spaces, thus extending the result  of Theorem 1.4 in \cite{OS}. (See also Theorem 2.1 in \cite{KMOS} and Theorem 41.4.1 in \cite{KM}.)

Our result also strengthens a result of Baldwin and Etnyre \cite{BE}.  They give a set of examples showing that when an open book decomposition has multiple binding components, no fixed lower bounds on fractional Dehn twisting can guarantee weak symplectic fillabilty.  This can now be viewed as a non-existence theorem for taut, oriented $C^{1,0}$-foliations satisfying our transitivity condition.

   The techniques introduced in this paper are used by the authors in \cite{KazezRoberts} to prove general $C^0$ approximation theorems; in particular, in Corollary~2.9, the condition that $\mathcal F$ be $\Phi$-bracketed can be replaced by the condition that $\mathcal F$ be taut, with continuous tangent plane field.
Closely related $C^0$ approximation results have also been obtained independently by J Bowden \cite{Bowden} using different techniques.

We thank Larry Conlon, John Etnyre, and Ko Honda for many helpful conversations.   We would also like to thank the referee for several helpful suggestions and corrections.

\section{An overview}

The transition from taut foliations to tight contact structures involves two auxiliary structures, volume preserving flows and symplectic topology.  We summarize the results we need from each field as follows.

\begin{thm} (Theorem~II.20, \cite{Sullivan}; see also Theorem A1, \cite{Hass})
Suppose $\mathcal F$ is a taut codimension-1 $C^{1,0}$-foliation of a smooth closed Riemannian 3-manifold $M$. Then there is a volume-preserving smooth flow $\Phi$ transverse to $\mathcal F$.
\end{thm}

For clarity, we break  the theorem  found in \cite{ET} into two statements:

\begin{thm} \label{weaklysymplectic} (Corollaries~3.2.2,  3.2.4 and 3.2.8, \cite{ET}; see also  Theorem~41.3.2, \cite{KM})  Let $M$ be a smooth closed Riemannian 3-manifold with a volume preserving flow $\Phi$.  Suppose there exist a smooth  positive  contact structure $\xi_+$ and a smooth negative contact structure $\xi_-$, both of which are transverse to $\Phi$.  Then each of $\xi_{\pm}$ is weakly symplectically fillable and universally tight.  Moreover,  if $\xi$ is any smooth (positive or negative) contact structure  transverse to $\Phi$, then $\xi$ is weakly symplectically fillable and universally tight.
\end{thm}

\noindent {\bf Remark:} The statement of Theorem~\ref{weaklysymplectic} is meant to emphasize two things.  First, given a smooth positive (respectively, negative) contact structure transverse to a volume preserving flow,  it is sufficient to produce a negative (respectively, positive) contact structure also transverse to the flow to conclude both are weakly symplectically fillable and universally tight.  Next, once such $\xi_+$ and $\xi_-$ are shown to exist, any contact structure $\xi$ transverse to $\Phi$ is necessarily weakly symplectically fillable and universally tight.

\begin{thm} \label{weaklysymplectic2} (Corollary~3.2.5, \cite{ET}; see also  Theorem~41.3.2, \cite{KM}) Let $M$ be a smooth closed oriented 3-manifold which contains a taut, oriented $C^2$-foliation $\mathcal F$. There exist a smooth positive contact structure $\xi_+$ and a smooth negative contact structure $\xi_-$, both $C^0$-close to $\mathcal F$.
\end{thm}

After giving background definitions and facts about foliations in \S\-~\ref{foliation basics}, we describe in \S\-\ref{transitive flow boxes} how flow boxes can be organized with an eye towards spreading an initial contact structure throughout the ambient manifold. Contact structures are propagated from one flow box to the next via a collection of local extension theorems described in \S\-\ref{basicsection}.  This leads to an inductive construction of the desired contact structure in \S\-\ref{theconstructionsection}.  Throughout, it is helpful to keep in mind the following

\medskip
\noindent{\bf Guiding Principle:} Constructions must be kept transverse to the flow. Moreover, when constructing a positive contact structure, the slope of the characteristic foliation of a partially constructed confoliation must be greater than or equal the slope of the intersection of the given foliation of $\mathcal F$ and vertical boundary of our flow boxes, with equality allowed only where $\mathcal F$ is smooth. When constructing a negative contact structure, the slope  inequality is reversed. \medskip

To explain this principle more formally, suppose that a closed oriented 3-manifold $M$ is expressed as a union of smooth submanifolds $V$ and $W$, possibly with corners, with $\partial V=\partial W$.  Suppose, moreover, that $W$ admits a codimension-1 foliation $\mathcal F_W$.  These submanifolds will be chosen so that their common boundary decomposes into {\sl horizontal} and {\sl vertical} portions, that is, portions tangent and transverse, respectively, to $\mathcal F_W$.  If a confoliation $\xi_V$ has been constructed on $V$ so that it is tangent to the horizontal portion of $\partial V$, transverse to the vertical portion, and contact on certain prescribed portions of $V$, then we call $V$ a {\sl contact zone}.

The Guiding Principle is a statement that the confoliation $\xi_V$ on $V$ must {\sl dominate} the foliation $\mathcal F_W$ along the vertical boundary (see Definition~\ref{compatible}) both for an initial choice of $V$, and also for subsequent choices as $V$ is expanded to all of $M$, and $W$ is shrunk correspondingly.

To expand a contact zone $V$ to the entire manifold $M$, we use the following structure.

\begin{definition} A foliation $\mathcal F_W$ is {\sl $V$-transitive} if every point in $W$ can be connected by a path in a leaf of $\mathcal F$ to a point of $V$.
\end{definition}

We will see in Theorem~\ref{main1} that the following structure is very useful.

\begin{definition} \label{tridecomposition} A closed 3-manifold $M$ admits a {\sl positive (respectively, negative) $(\xi_V,\mathcal F_W,\Phi)$ decomposition} if $M$ can be decomposed as  a union 
$$M=V\cup W,$$
where  the horizontal portion of $\partial W$ is tangent to $\mathcal F_W$ and the vertical portion of $\partial V$ is tangent to $\Phi$, and
\begin{enumerate}
\item $\mathcal F_W$ is a $V$-transitive oriented foliation of $W$,

\item $\xi_V$ is a smooth contact structure defined on $V$ which positively (respectively, negatively) dominates $\mathcal F_W$,

\item for some choice of Riemannian metric, $M$ admits a volume preserving flow $\Phi$ transverse to both $\xi_V$ and $\mathcal F_W$.

\end{enumerate}
\end{definition}

Note that the existence of a $(\xi_V,\mathcal F_W,\Phi)$ decomposition does not require the existence of a codimension-1 foliation defined on   all of $M$.

\medskip

\noindent{\bf Theorem~\ref{main1}.}
{\it If $M$ admits a positive $(\xi_V,\mathcal F_W,\Phi)$ decomposition, then $M$ admits a smooth positive contact structure  which agrees with $\xi_V$ on $V$ and is  $\Phi$-close to $\mathcal F_W$ on $W$. The analogous result holds if $M$ admits a negative $(\xi_{V'},\mathcal F_{W'},\Phi)$ decomposition. If $M$ admits both a positive $(\xi_V,\mathcal F_W,\Phi)$ decomposition and a negative $(\xi_{V'},\mathcal F_{W'},\Phi)$ decomposition, then these   contact structures $(M,\xi_+)$ and $(-M,\xi_-)$ are weakly symplectically fillable and universally tight.}
\medskip

\begin{question} If a closed oriented 3-manifold admits both a positive $(\xi_V,\mathcal F_W,\Phi)$ and a negative $(\xi_{V'},\mathcal F_{W'},\Phi)$ decomposition, does it contain a taut oriented foliation transverse to $\Phi$?
\end{question}

\begin{definition} Given a splitting $M=V\cup W$ and a flow $\Phi$, an oriented, codimension-1 foliation $\mathcal F$  is {\sl compatible} with $(V, W, \Phi)$ if $\mathcal F$ is transverse to $\Phi$, and the common boundary $\partial V=\partial W$ decomposes into subsurfaces which are either horizontal or vertical with respect to $\mathcal F$.
\end{definition}

\begin{definition} An oriented codimension-1 foliation $\mathcal F$ of a 3-manifold $M$ is {\sl bracketed} if, for some volume preserving flow $\Phi$,
\begin{enumerate}

\item $\mathcal F$ is compatible with some $(V,W,\Phi)$ decomposition of $M$ for which there exist $\mathcal F_W$ and $\xi_V$ such that $(\xi_V,\mathcal F_W,\Phi)$ is a positive decomposition, and

\item $\mathcal F$ is compatible with some $(V',W',\Phi)$ decomposition of $M$ for which there exist $\mathcal F_{W'}$ and $\xi_{V'}$ such that $(\xi_{V'},\mathcal F_{W'},\Phi)$ is a negative decomposition.
\end{enumerate}
When we wish to specify the flow $\Phi$, $\mathcal F$ is called {\sl $\Phi$-bracketed}.
\end{definition}

\begin{cor}
Let $\mathcal F$ be an oriented codimension-1 foliation of a 3-manifold $M$ which is $\Phi$-bracketed. Then there exist a smooth positive contact structure $\xi_+$ and a smooth negative contact structure $\xi_-$, both $\Phi$-close to $\mathcal F$.  These contact structures $(M,\xi_+)$ and $(-M,\xi_-)$ are weakly symplectically fillable and universally tight.
\end{cor}

Sometimes $\mathcal F_W$ and $\mathcal F_{W'}$ are obtained by restricting $\mathcal F$ to $W$ and $W'$ respectively, and sometimes they are not. Very roughly speaking, when $W=W'$, we think of the restriction of $\mathcal F$ to $W$ as being {\sl bracketed by} $\mathcal F_{W'}$ and $\mathcal F_W$  as a generalization of the situation in which the slope of $\mathcal F$ along boundary components of $W$ lies between the corresponding boundary slopes of $\mathcal F_{W'}$ and $\mathcal F_W$.

As noted in Corollary~\ref{smoothisbracketed}, all taut, oriented $C^2$-foliations apart from the product foliation $S^1\times S^2$ are bracketed.  In Example~\ref{s1timess2}, we show that the product foliation of $S^1\times S^2$ is not bracketed. In this paper, we show that many taut, oriented $C^{1,0}$-foliations are bracketed. For each bracketed foliation considered in this paper, it is possible to choose $V=V'$ and $W=W'$.

\begin{con}\label{bracketconj}
Let $\mathcal F$ be a taut oriented $C^{1,0}$-foliation of a closed oriented 3-manifold $M\ne S^1\times S^2$. Then $\mathcal F$ is bracketed.
\end{con}

We have the following two closely related questions.

\begin{question} Suppose $\mathcal F$ is a taut, oriented foliation with no torus leaf. Is $\mathcal F$ $C^0$-close to a taut, oriented smooth foliation? 
\end{question}

\begin{question} Suppose $\mathcal F$ is a taut, oriented foliation with no torus leaf.  Is $\mathcal F$ $\Phi$-close to a taut, oriented smooth foliation for some volume preserving flow $\Phi$?
\end{question}

\medskip

Establishing the initial contact zone is of fundamental importance. In the context of $C^{1,0}$-foliation theory we introduce, in \S\-\ref{attracting holonomy}, the notion of {\sl holonomy with a contracting interval} and define what we mean by {\sl attracting neighborhood}.  This is significantly weaker than the more familiar notion of linear attracting holonomy, yet it suffices to build an initial contact zone. The precise definition appears as Definition~\ref{contracting}. As a corollary to Theorem~\ref{main1}, we obtain:

\medskip
\noindent{\bf Theorem~\ref{attracting}.}
{ Let $\mathcal F$ be a taut $C^{1,0}$-foliation transverse to a flow $\Phi$.  If $V$ is a disjoint union of attracting neighborhoods, and $\mathcal F$ is $V$-transitive, then $\mathcal F$ is bracketed and hence can be $\Phi$-approximated by a pair of  contact structures  $\xi_{\pm}$, one positive and one negative.  These contact structures $(M,\xi_+)$ and $(-M,\xi_-)$ are weakly symplectically fillable and universally tight.
}
\medskip

When working with $C^{1,0}$-foliations, it can be difficult, or even impossible, to establish the existence of sufficient  attracting holonomy. Therefore, we introduce a different way of creating an initial contact zone. Roughly speaking, instead of looking for loops tangent to the foliation and satisfying a nice property, we look for loops \textit{transverse} to the foliation and satisfying a nice property.  We make this precise in Definition~ \ref{transtranv}, where we define  {\sl $L$-bracketed foliation}. As a corollary to Theorem~\ref{main1}, we obtain:
 
\medskip
\noindent{\bf Theorem~\ref{transitivemain}.}
{\it Suppose $\mathcal F$ is a taut oriented codimension-1 foliation in $M$, and that $\mathcal F$ is $L$-bracketed for some link $L$. Then $\mathcal F$ is bracketed and hence can be $\Phi$-approximated by a pair of smooth contact structures $\xi_{\pm}$, one positive and one negative. These contact structures are necessarily weakly symplectically fillable and  universally tight.
}
\medskip

In \S~\ref{OBresults}, we consider the important special case that $\mathcal F$ is transverse to a flow $\Phi$ that has been obtained by removing a fibred link $L$ and doing a Dehn filling of a volume preserving suspension flow. In this case, $L$ forms the binding of an open book decomposition $(S,h)$ of $M$ and the contact structure $\xi_{(S,h)}$ compatible with $(S,h)$ is $\Phi$-close to $\mathcal F$.  

In \cite{HKM2}, Honda, Kazez and Mati\'c introduced the use of  foliations $\Phi$-close to $\xi_{(S,h)}$ as a way of establishing universal tightness of $\xi_{(S,h)}$. In particular, they appealed to $C^{1,0}$-foliations constructed in \cite{R1,R2} to claim that $\xi_{(S,h)}$ is universally tight whenever the binding of $(S,h)$ is connected and the fractional Dehn twist coefficient at least one. Although the foliations constructed in \cite{R1,R2} are not smooth, and therefore the proof in \cite{HKM2} contained a gap, they are $L$-bracketed, and hence 
Theorem~\ref{transitivemain} reveals that the conclusions of  \cite{HKM2}  are correct. 

In \S\-\ref{Open book} we also  include some background material relating language arising in the theory of open books with language arising in the theory of  foliations. In particular, we give a translation between coordinates used in each subject together with a summary of our results related to open book decompositions.  

To make the paper more self-contained there is an appendix containing an overview of the relationship between volume preserving flows and closed dominating 2-forms, and giving some standard definitions from symplectic topology. Most of this material is present either implicitly or explicitly in \cite{ET}.

We close this section with an application of Theorem~\ref{main1} to the study of L-spaces.

\begin{definition} (Definition 1.1, \cite{OS2})
A closed three-manifold $Y$  is called an {\sl L-space} if $H_1(Y;\mathbb Q)=0$ and $\widehat{HF}(Y)$ is a free abelian group of rank $|H_1(Y;\mathbb Z)|$.
\end{definition}

\begin{thm} (Theorem 1.4, \cite{OS})
An L-space has no symplectic semi-filling with disconnected boundary; and all its symplectic fillings have $b_2
(W)=0$. In particular, $Y$ admits no taut smooth foliation.
\end{thm}
 
In other words, Ozsv\'ath and Szab\'o show that if $Y$ is an L-space then there is no symplectic manifold $(X,\omega)$ with weakly convex boundary such that $|\partial X|>1$ and $Y$ is one of the boundary components. So an L-space cannot contain a pair of $\Phi$-close contact structures, $\xi_+$ positive and $\xi_-$ negative, where $\Phi$ is a volume preserving flow.
(For details, see the Appendix.) Theorem~\ref{main1} thus implies the following.

\begin{cor} An L-space $Y$ admits no bracketed foliation.
\end{cor} 

In particular, the foliations constructed in \cite{KRo, LR, R1, R2} never exist in an L-space.

\section{Foliation basics}\label{foliation basics}

\begin{definition}\label{folndefn}
Let $M$ be a smooth closed 3-manifold, and let $k$ be a non-negative integer or infinity.  A {\sl $C^k$ codimension-1 foliation} $\mathcal F$ of (or in) $M$ is a union of disjoint  connected  surfaces $L_i$, called the {\sl leaves} of $\mathcal F$, such that:
\begin{enumerate}
\item $\cup_i L_i = M$, and
\item  there exists a $C^k$ atlas $\mathcal A$ on $M$ which contains all $C^{\infty}$ charts and with respect to which $\mathcal F$ satisfies the following local product structure: 
\begin{itemize}
\item[] for every $p\in M$, there exists a coordinate chart $(U,(x,y,z))$ in $\mathcal A$ about $p$ such that $U\approx \mathbb R^3$ and the restriction of $\mathcal F$ to $U$ is the union of
planes given by $z = $ constant. 
\end{itemize}
\end{enumerate}

  If the tangent plane field $T\mathcal F$ exists and is continuous, $\mathcal F$ is called $C^{1,0}$. If $\mathcal F$ is $C^{1,0}$ and all leaves are smoothly immersed, then $\mathcal F$ is called $C^{\infty,0}$. (Compare with Definition~1.2.22 of \cite{CC} and see Definition~2.2 of \cite{KazezRoberts}.)
\end{definition}

 When $T\mathcal F$ is continuous,  there is a continuous, and hence  a smooth, 1-dimensional foliation transverse to $\mathcal F$. Moreover, it follows by Theorem~\ref{calegarikr} that if $\mathcal F$  is $C^{1,0}$, then it can be $C^0$ approximated by a $C^{\infty,0}$ foliation.  Thus in seeking to approximate a $C^{1,0}$-foliation, no generality is lost in assuming the leaves are smoothly immersed.

A frequently used technique for constructing foliations is to start with a branched surface embedded in $M$ that has product complementary regions.  Since the embedding may be smoothed, a foliation resulting from thickening the branched surface and extending across the complementary regions can be constructed to be $C^0$.

Definition~\ref{folndefn} extends in an obvious way to define a codimension-1 foliation on a compact oriented smooth 3-manifold with non-empty boundary, where we insist that for each torus boundary component $T$, either $T$ is a leaf of $\mathcal F$, or $\mathcal F$ is everywhere transverse to $T$, and that any non-torus boundary component is a leaf of $\mathcal F$.

Recall that a  smooth structure with corners on a topological 3-manifold $M$ with nonempty boundary is a maximal collection of smoothly compatible charts with corners whose domains cover $M$, where a chart  with corners is an open set diffeomorphic to one of $\mathbb R^3$, $\{(x,y,z)\}| z\ge 0\}$, $\{(x,y,z)\}| y,z\ge 0\}$, or $\{(x,y,z)\}| x,y, z\ge 0\}$. Notice that the boundary of a manifold with corners naturally admits a stratification as a disjoint union of 0-, 1-, and 2-dimensional manifolds. The  0- and 1-manifolds of this stratification are referred to as the corners of $M$.

Definition~\ref{folndefn}  extends further in an obvious way to define a codimension-1 foliation on a compact smooth 3-manifold $M$ with corners, where we insist that $\partial M$ can be written as a union of two compact piecewise linear surfaces $\partial_v M$ and $\partial_h M$, where the intersection $\partial_v M\cap\partial_h M$ is a union of corners of $M$, the components of $\partial_h M$ are contained in leaves of $\mathcal F$, and $\partial_v M$ is everywhere transverse to $\mathcal F$.

\begin{definition} A {\sl flow} is an oriented 1-dimensional foliation of $M$; namely, a decomposition $\Phi$ of a smooth compact 3-manifold  $M$ into a disjoint  union of connected  1-manifolds, called the {\sl flow curves} of $\Phi$, such that  there exists a $C^k$ atlas $\mathcal A$ on $M$ which contains all $C^{\infty}$ charts and with respect to which $\Phi$ satisfies the following local product structure: 
\begin{itemize}
\item[] for every $p\in M$, there exists a coordinate chart $(U,(x,y,z))$ in $\mathcal A$ about $p$ such that $U\approx \mathbb R^3$, and the restriction of $\Phi$ to $U$ is the union of
lines given by $(x,y) = $ constant. 
\end{itemize}
When $M$ has boundary a disjoint union of tori, we insist that for each torus boundary component $T$, either $\Phi$ is everywhere tangent to $T$  or $\Phi$ is everywhere transverse to $T$.
When $M$ is smooth with corners, we insist that $\partial M$ can be written as a union of two compact piecewise linear surfaces $\partial_v M$ and $\partial_h M$, where the intersection $\partial_v M\cap\partial_h M$ 
is a union of corners of $M$,  $\partial_h M$ is everywhere transverse to $\Phi$, and $\partial_v M$ is everywhere tangent to $\Phi$.
\end{definition}

Flows and oriented codimension-1 foliations coexist in interesting ways. A good overview can be found in \cite{CC}.  In particular, given an oriented   $C^{1,0}$-foliation $\mathcal F$  of an oriented 3-manifold $M$, possibly with non-empty boundary and possibly with corners, there is a $C^{\infty}$ flow everywhere transverse to $\mathcal F$. From this we have the following:

 \begin{thm} \cite{calegari,KR4}\label{calegarikr}
Suppose $\mathcal F$ is a $C^{1,0}$ foliation in $M$.  Then there is an isotopy of $M$ taking $\mathcal F$ to a $C^{\infty, 0}$ foliation $\mathcal G$ with smoothly immersed leaves which is $C^0$ close to $\mathcal F$.  If $\Phi$ is a smooth flow transverse to $\mathcal F$, the isotopy may be taken to map each flow line of $\Phi$ to itself.  
\end{thm}

When the oriented foliation $\mathcal F$ is taut, and $M$ is Riemannian, the smooth transverse flow can be chosen to be volume preserving.

\begin{thm} (Theorem~II.20, \cite{Sullivan}; see also Theorem A1, \cite{Hass})\label{volpreserve} Let $\mathcal F$ be a codimension-1, taut $C^{1,0}$-foliation of a closed smooth Riemannian 3-manifold $M$. Then there is a volume-preserving smooth flow everywhere transverse to $\mathcal F$. Equivalently, there is a smooth closed 2-form dominating $\mathcal F$.
\end{thm}

Given a 3-manifold $M$ containing a taut $C^{1,0}$-foliation, one can ask whether there is a closely related $C^{\infty}$-foliation. Interpreting `closely related'  to mean any of  $C^0$-$\epsilon$-close for some fixed $\epsilon>0$, $\Phi$-close, or topologically conjugate results in questions for which the answers are very little understood. There are certainly 3-manifolds which contain Reebless $C^{1,0}$-foliations but not Reebless $C^{\infty}$-foliations (Theorem D, \cite{BNR}).  The existence of a taut sutured manifold hierarchy guarantees the existence of two types of foliation, one   $C^0$ and finite depth and the other $C^{\infty}$ (\cite{G1,G2,G3}).  Fixing a Riemannian metric and some $\epsilon>0$, these two types of foliation are not necessarily $C^0$-$\epsilon$-close. However, since they are carried by a common transversely oriented branched surface, they are $\Phi$-close.
 
We will take advantage of the fact that  it is always possible {\sl locally} to $C^0$-approximate $(\mathcal F,T\mathcal F)$ by   $(\tilde{\mathcal F}, T\tilde{\mathcal F})$, for some locally defined smooth foliation  $\tilde{\mathcal F}$.

\begin{prop} \label{smoothapprox}  Let $D$ be a smooth disk with corners  and let $\mathcal F$ be a $C^{1,0}$-foliation of $D^2\times [0,1]$  which is positively transverse to the smooth 1-dimensional foliation by intervals $\{ (x,y)\}\times [0,1],\,\, (x,y)\in D^2$. Given any $\epsilon>0$, there is a smooth foliation $\tilde{\mathcal F}$ which is  positively transverse to the smooth 1-dimensional foliation by intervals $\{ (x,y)\}\times [0,1],\,\, (x,y)\in D^2$, and satisfies $(\tilde{\mathcal F},T\tilde{\mathcal F})$ is  $C^0$ $\epsilon$-close to $(\mathcal F, T\mathcal F)$. Moreover, if $\mathcal F$ is smooth on some compact $\mathcal F$-saturated subset, then we may choose $\tilde{\mathcal F}$ to equal $\mathcal F$ on this subset.
\end{prop}

\begin{proof}  By identifying $D$ with a subset of the plane, a point $p$ in a leaf of $\mathcal F$ determines both a point in $\mathbb R^3$, and by choosing a unit vector perpendicular to the tangent plane of the leaf, a point $\bf{u}_p$ in $T\mathbb R^3$.  The standard metric on $T\mathbb R^3 = \mathbb R^6$ is used to measure the distance between two leaves of $\mathcal F$ as follows.

By  Theorem~\ref{calegarikr} , we may assume that the leaves of  $\mathcal F$ are given by  the graphs of $z=f_{\theta}(x,y)$, for some continuous family of smooth functions $f_{\theta} : D\to [0,1]$, $0\le\theta\le 1$.
For any two such leaves, $L_1$ and $L_2$ say, given by $z=f_{\theta_1}(x,y)$ and $z=f_{\theta_2}(x,y)$ respectively, define the distance between them to be the maximum distance, computed in $T\mathbb R^3$, between  $(x,y,f_{\theta_1}(x,y),{\bf u}_{(x,y,f_{\theta_1}(x,y))})$ and $(x,y,f_{\theta_2}(x,y),{\bf u}_{(x,y,f_{\theta_2}(x,y))})$ for $(x,y) \in D$.  Since $D$ is compact, uniform continuity guarantees  that $d$ is continuous and hence a metric on the leaf space of $\mathcal F$.

For any $\theta \in [0,1]$, let $U_{\theta}$ denote the subset of $D\times [0,1]$ which is the union of all graphs $z=f_{\theta}(x,y)$ which are of d-distance strictly less than $\epsilon/2$ from the leaf $z = f_{\theta}(x,y)$. Since $U_{\theta}$ is the pullback of an $\epsilon/2$ $d$-neighborhood in the leaf space, $U_{\theta}$ is open in $D\times[0,1]$.

 Pick a finite cover of $D\times[0,1]$ by $U_{\theta_0},U_{\theta_1},\cdots,U_{\theta_r}$ for some $r\ge 0$ and  $0=\theta_0<\theta_1<\cdots <\theta_r = 1$. Now let $\tilde{\mathcal F}$ be the foliation of $D\times[0,1]$ which includes the leaves given by the graphs of $f_{\theta_i}$ and, for each $i, 0\le i\le r-1$, the leaves given by  the graphs of a damped straight line homotopy between $f_{\theta_i}$ and $f_{\theta_{i+1}}$.  Thus if $g$ is a smooth homeomorphism of $[0,1]$ with derivatives at $0$ and $1$ vanishing to infinite order, the leaves of $\tilde{\mathcal F}$ are $z=(1-g(t))f_{\theta_i}(x,y) + g(t)f_{\theta_{i+1}}(x,y)$, $0\le t\le 1$, on the subset of $D\times [0,1]$ bounded by the graphs of $f_{\theta_i}$ and $f_{\theta_{i+1}}$.   By construction, $\tilde{\mathcal F}$ is smooth.

Moreover, $(\tilde{\mathcal F},T\tilde{\mathcal F})$ and $(\mathcal F,T\mathcal F)$ are $\epsilon$-close. To see this, recall that a normal vector to a graph $z=f(x,y)$ is given by ${\bf n}_f = \langle -f_x, -f_y, 1\rangle$, and a straight-line homotopy between $f_{\theta_1}$ and $f_{\theta_2}$ induces a straight-line homotopy between ${\bf n}_{f_{\theta_1}}$ and ${\bf n}_{f_{\theta_2}}$. Normalizing this straight-line homotopy of normal vectors gives a geodesic on the unit sphere joining ${\bf u}_{f_{\theta_1}}$ and ${\bf u}_{f_{\theta_2}}$. Since the leaves given by $z = f_{\theta_i}(x,y)$ and $z = f_{\theta_{i+1}}(x,y)$ are of $d$-distance at most $\epsilon/2$, it follows immediately from the triangle inequality that the leaves given by $z=f_{(1-g(t))\theta_i + g(t)\theta_{i+1}}(x,y)$ and $z=(1-g(t))f_{\theta_i}(x,y) + g(t)f_{\theta_{i+1}}(x,y)$ are of $d$-distance strictly less than $\epsilon$. So$(\tilde{\mathcal F},T\tilde{\mathcal F})$ and $(\mathcal F,T\mathcal F)$ are $\epsilon$-close.

If $\mathcal F$ is smooth on some compact $\mathcal F$-saturated subset $A$ of $D\times [0,1]$, each component of $A$ is bounded by graphs of the form $f_{\theta}$. By compactness of $A$, $\partial A$ contains only finitely many such $f_{\theta}$. For each $z=f_{\theta}$ in $\partial A$, include $\theta$ in the list $\theta_0,\theta_1,\cdots,\theta_r$ and modify $\mathcal F$ only on the  complement of $A$.
\end{proof}

Next we recall Operation 2.1.1 of  \cite{Ga}. Let $L_1,\dots, L_m$ be distinct leaves of 
a $C^{1,0}$-foliation $\mathcal F$. Modify $\mathcal F$ by thickening each of the leaves $L_j$. Thus, each $L_j$ is blown up to an $I$-bundle $L_j\times [-1,1]$.  Let $\mathcal F'$ denote the resulting foliation. We highlight the following observation.

\begin{lemma} \label{blowup}  The leaves $L_1, \dots, L_m$ may be thickened so that  the foliation $\mathcal F'$ is $C^0$ and the restriction of $\mathcal F'$ to  $$L_j\times (-1,1)\subset M$$ is a smooth foliation for each $j$.
\end{lemma}

\section{Transitive flow box decompositions}\label{transitive flow boxes}

\begin{definition} \label{flowboxdefn}

A {\sl flow box}, $F$, for a $C^{1,0}$-foliation $\mathcal F$ and smooth transverse flow $\Phi$, is a smooth chart with corners that is of the form $D\times I$, where $D$ is a polygon (a disk with at least three corners), $\Phi$ intersects $F$ in the arcs $\{(x,y)\}\times I$, and $\mathcal F$  intersects $F$ in disks which are everywhere transverse to  $\Phi$ and hence can be thought of as graphs  over $D$. In particular, $D\times \partial I$ lies in leaves of $\mathcal F$, each component of $\mathcal F\cap F$ is a smoothly embedded disk, and these disks vary continuously in the $I$ direction.

The {\sl vertical boundary} of $F$, denoted $\partial_v F$, is $\partial D \times I$.  The {\sl horizontal boundary} of $F$ is $D \times \partial I$ and is denoted $\partial_h F$.  An arc in $M$ is {\sl vertical} if it is a subarc of a flow line and {\sl horizontal} if it contained in a leaf of $\mathcal F$. 

It is often useful to view the  disk $D$ as a 2-cell with $\partial D$ the cell complex obtained by letting the vertices correspond exactly to the  corners of $D$. Similarly, it is useful to view the flow box $F$ as a 3-cell  possessing the product cell complex structure of $D\times I$. Then
the horizontal boundary $\partial_h F$ is  a union of two (horizontal) 2-cells and the vertical boundary $\partial_v F$ is a union of $c$ (vertical) 2-cells, where $c$ is the number of corners of $D$.  

A subset $R$ of $F$ is called a {\sl vertical rectangle} if it has the form $\alpha\times [a,b]$, where $0\le a<b\le 1$ and $\alpha$  is either a 1-cell of $\partial D$ or else  a properly embedded arc in $D$ connecting distinct vertices of $D$. A subset $e$ of $F$ is called an {\sl  edge} if it is a compact interval contained in a 1-cell of $F$.

Given a vector $\vec w$ tangent to $\partial_v F$, we choose a {\sl slope} convention such that the leaves of $\mathcal F \cap \partial_v F$ have slope $0$, the flow lines have slope $\infty$, and the sign of the slope of $\vec w$ is computed as viewed from outside of $F$.
\end{definition}

\begin{definition}  Given a codimension-1 leafwise smooth foliation $\mathcal F$ and transverse smooth flow $\Phi$, let $V$ be a compact codimension-0 sub-manifold of $M$, with $\partial V = \partial_v V \cup\partial_h V$, where $\partial_v V$ is a union of flow arcs or circles, and $\partial_hV$ is a union of subsurfaces of leaves of $\mathcal F$. In the case that $\partial V = \partial_vV$, $\mathcal F$ and $\Phi$ need only be defined on the complement of $V$.  A {\sl flow box decomposition} of $M$ {\sl rel} $V$ is a decomposition of $M\setminus \text{int} V$ as a finite  union $M = V\cup (\cup_iF _i)$ where

\begin{enumerate}
\item Each $F _i$ is a standard flow box for $\mathcal F$.

\item If $i \neq j$, the interiors of $F _i$ and $F _j$ are disjoint.

\item If $F _i$ and $F _j$ are different flow boxes, then their intersection is connected and either empty, a  0-cell, an edge, a  vertical rectangle, or a subdisk of $\partial_h F_i \cap \partial_h F_j$.  

\end{enumerate}
\end{definition}

\begin{definition} \label{transitive} We call a flow box decomposition $M= V\cup F_1 \cup \dots \cup F_n$ {\it transitive} if $V_0=V$, $V_i = V_{i-1} \cup F_i$, and for $i=1,\dots, n$,

\begin{enumerate}
\item each 2-cell  of $\partial_v F_i$ has interior disjoint from $\partial_h F_j$ for all $j<i$,
\item $V_{i-1}\cap F_i$ is a union of horizontal subsurfaces and vertical 2-cells of $F_i$, together possibly with some 0- and 1-cells, 
and 

\item $V_{i-1}\cap F_i$ contains a vertical 2-cell of $F_i$.
\end{enumerate}

\end{definition}

\begin{prop}\label{transitiveflowbox}   If $M$ is $V$-transitive, then there is a transitive flow box decomposition of $M$ rel $V$.
\end{prop}

\begin{proof} 
Since $M$ is $V$-taut, for each point $x\in M\backslash V$, there exists $\gamma$ an embedded arc in the leaf containing $x$ that connects $x$ to $V$.  By taking a regular neighborhood of $\gamma$ in its leaf and flowing along it, create a flow box $F$ with the property that it has one vertical 2-cell contained in $\partial_v V$.  The point $x$ may or may not be in the interior of $F$.
Using compactness of $M$, pick a finite collection of flow boxes, of the sort just described, $F_1, F_2,\dots, F_r$ that cover $M\backslash V$. Assume no proper subset of the $F_i, 1\le i\le r,$ cover.  

Next, let $L_1, L_2, \dots, L_m$ be the collection of leaves of $\mathcal F$ that contain the horizontal boundaries of all $F_i$.  

We proceed by induction to show that $V\cup F_1\cup\dots\cup F_i$ admits a flow box decomposition with respect to $V$ for every $i, 1\le i\le r$.  Certainly, $V\cup F_1$ does. So suppose that $V\cup F_1\cup\dots\cup F_{i-1}$ admits a flow box decomposition with respect to $V$. After renaming and reindexing as necessary, we  assume that this flow box decomposition is given by $V\cup F_1\cup\dots\cup F_{i-1}$.  We show that $V_i = V\cup F_1\cup\dots\cup F_{i-1}\cup F_i$ also admits a flow box decomposition with respect to $V$.

Begin by slightly increasing the size of  $F_i$ in $M\setminus \text{int} V$, as necessary, so that $F_i$ is still a flow box and, for all $j<i$, $\partial_v F_j$ and $\partial_v F_i$ are transverse  away from $V$, where they may overlap tangentially. Notice that this ensures that $ V\cup F_1\cup\dots\cup F_{i-1}\cup F_i$ is a codimension-0 submanifold with corners and  piecewise vertical and horizontal boundary. Also, cut $F_i$ open along those (horizontal disk) components of $(\cup_i L_i)\cap F_i$ which have non-empty intersection with $\partial_h F_j$ for some $j<i$.  Denote the resulting flow boxes by $F_i^1,\dots,F_i^s$; so $F_i = F_i^1\cup\dots\cup F_i^s$.

Consider $F_j\cap F_i^1$ for some $j<i$. Since $\partial_v F_j$ and $\partial_v F_i^1$ are transverse away from $V$, each component of $F_j\cap F_i^1$ is  a flow box. Consider any such component, $X$ say, from the point of view of the flow box  $F_i^1=D_i\times [c,d]$. Notice that  $X=D\times [c,d]$, where $D$ is a subdisk (with corners) of $D_i$, and $X\cap \partial_v F_i^1$ is a non-empty union  of vertical 2-cells. 

Now, for all $j<i$, remove $F_j\cap F_i^1$  from $F_i^1$. Taking the closure of the result we get  a union $G_1\cup\dots\cup G_b$ of flow boxes, where each  $G_k\cap   V_{i-1}$ is a union of horizontal subsurfaces and vertical 2-cells of $G_k$, together possibly with some 0- and 1-cells, and contains a vertical 2-cell of $G_k$.  Notice that $G_k\cap G_l=\emptyset$ if $k\ne l$ and, by subdividing each $G_k$ along finitely many vertical rectangles as necessary, we may assume $G_k\cap F_j$ is connected for all $j<i$.  The resulting union
$$V\cup F_1\cup\dots\cup F_{i-1}\cup G_1\cup\dots\cup G_b$$ is then a transitive flow box decomposition of $V\cup F_1\cup\dots\cup F_{i-1}\cup F_i^1$ with respect to $V$.

Repeat this process for each $a, 2\le a\le s$, to obtain a transitive flow box decomposition of $V\cup F_1\cup\dots\cup F_i$ with respect to $V$.
\end{proof}

\section{Basic extension results}\label{basicsection}

In this section, we collect together an assortment of confoliation extension results important for  the inductive construction to be described in Section~\ref{theconstructionsection}. For the most part, it will be possible to restrict attention to flow boxes diffeomorphic to one of the flow boxes $F$, $G$ or $H$, where $F$, $G$ and $H$ are defined as follows.

Let $F$ denote the flow box given by 
$$F=\{|x| \le 1, 0 \le y \le 1, |z| \le 1\}.$$
Let $G$ denote the flow box given by $\Delta \times [0,1]$, where $\Delta$ is the region in the $xy$-plane bounded by the triangle with vertices
$$ (-3/2,1), (3/2,1) \mbox{ and } (0,-1/2).$$ Let $\Delta^{(0)}$ denote the $0$-skeleton 
$$\Delta^{(0)} = \{(-3/2,1), (3/2,1), (0,-1/2)\}.$$
Let $H=F\cap G$, a flow box with hexagonal horizontal cross-section. Notice that $H$ is diffeomorphic to the complement in $G$ of an open neighborhood of the $1$-skeleton of 
$\partial_v G$.

We begin with the following elementary, and very useful, observations of Eliashberg-Thurston, \cite{ET}.

\begin{lemma} (Proposition~1.1.5, \cite{ET})\label{ETlemma} Let $\eta$ be a $C^k$-confoliation  with $k \ge 1$ and domain $F$ given by a 1-form $dz - a(x,y,z) dx$.  Then $$\dfrac {\partial a}{\partial y}(x,y,z) \ge 0$$ at all points of $F$, and $$\dfrac {\partial a}{\partial y}(x,y,z) > 0$$ where $\eta$ is contact. \qed
\end{lemma}

Eliashberg and Thurston use this lemma to approximate the confoliation by a contact structure with the following corollary.

\begin{cor} (Lemma~2.8.2, \cite{ET})\label{approximate} Let $\eta$ be a $C^k$-confoliation, with $k \ge 1$ and domain $F$, given by a 1-form $dz - a(x,y,z) dx$. If $\eta$ is contact in a neighborhood of $y=1$ in $F$, then $\eta$ can be approximated by a confoliation $\hat\eta$ which coincides with $\eta$ together with all of its derivatives along the boundary $\partial F$ and which is contact inside $F$.
\end{cor}

\begin{proof} It is enough to approximate $a(x,y,z)$ along each interval $\{x\} \times [0,1] \times \{z\}$ by a function $\hat{a}(x,y,z)$ that is strictly monotonic  for $(x,z)\in(-1,1) \times (-1,1)$ but is damped to agree smoothly with $a(x,y,z)$ on $\partial F$.
\end{proof}

\begin{definition}\label{dominate} If $\alpha$ and $\beta$ are families of curves transverse to $\partial/\partial z$ and contained in a vertical 2-cell $R$ of the vertical boundary of a flow box, we say $\alpha$ {\sl strictly dominates} $\beta$ along $A\subset R$ if at every $p\in A$, the slope of the tangent to $\alpha$ at $p$ is greater than the slope of the tangent to $\beta$ at $p$.  It must be specified if the comparison of slopes is made from inside or outside of the flow box. If $\alpha$ strictly dominates $A$, and $\alpha$ and $\beta$ are the characteristic foliations of 2-plane fields $\xi_1$ and $\xi_2$ respectively, we also say that $\xi_1$ strictly dominates $\xi_2$ along $R$. If $\xi_2 = T\mathcal F$ for some codimension-1 foliation $\mathcal F$, we also say that $\xi_1$ strictly dominates $\mathcal F$ along $R$.
\end{definition}

The statement of Lemma~\ref{ETlemma} raises the question of whether flow box coordinates can always be chosen so that the contact form can be written as $dz-a(x,y,z)dx$.  The next lemma points out that this is the case and gives a simple condition for a contact structure to  dominate in such coordinates.

\begin{lemma}
Let $U$ be a regular neighborhood in $F$ of the union of $x=\pm 1$ and $z=\pm 1$.  Let $\eta$ be a $C^k$-confoliation  with $k \ge 1$ defined in a neighborhood $V$ of $y=1$ in $F$ 
which is everywhere transverse to the vertical segments $(x,y) =$ constant, horizontal in $\overline{U}$, and contact on $V\setminus \overline{U}$. Then, after  smoothly reparametrizing  $F$ as necessary, we may assume that $\eta$ is given by a 1-form 
$$dz - a(x,y,z) dx$$ with 
\begin{enumerate}
\item  $a(x,y,z)=0$ on $V\cap \overline{U}$, and 
\item $\dfrac {\partial a}{\partial y}(x,y,z) > 0$ on $V\setminus \overline{U}$.
\end{enumerate}
Moreover, the characteristic foliation of $\eta$ along the complement of $\overline{U}$ in $y=1$ strictly dominates the horizontal foliation, when viewed from inside $F$, if and only if
$a(x,y,z)>0$ in $V'\setminus \overline{U}$, for some neighbourhood $V'\subset V$ of $y=0$. 
\end{lemma}

\begin{proof} Since $dz - a(x,y,z) dx$ vanishes on $\partial/\partial y$, it is enough to choose coordinates $x, y$ for leaves so that curves with constant $x$ coordinate are Legendrian.  At points where $\eta$ is transverse to the horizontal foliation, there is a unique Legendrian direction. At all other points, any direction is Legendrian.  The coordinate $y$ can be constructed by choosing a section of the Legendrian directions.
\end{proof}

\begin{figure}[htbp] %  figure placement: here, top, bottom, or page
   \centering
   \includegraphics[width=5in]{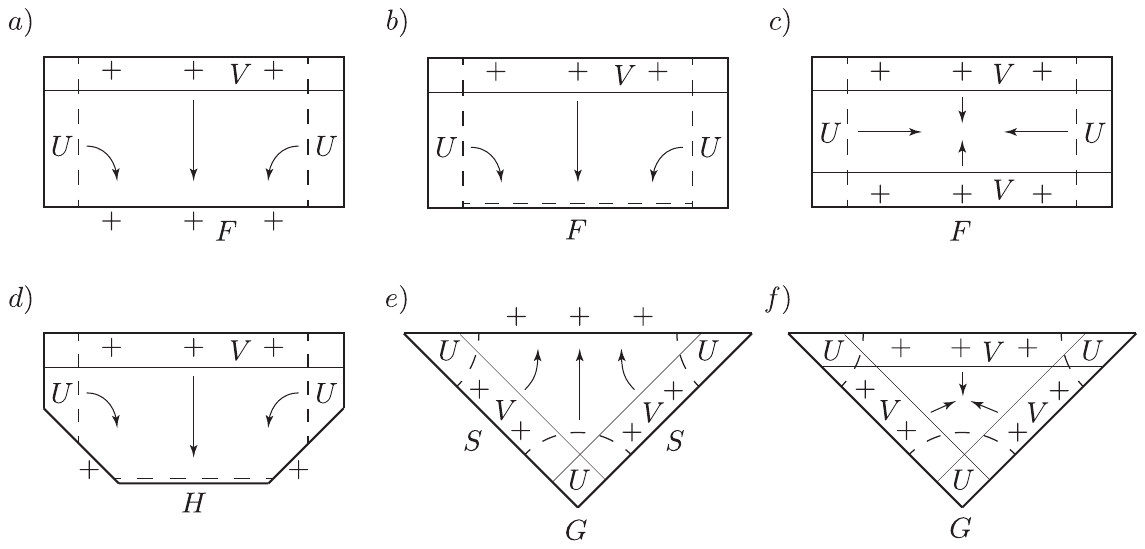} 
   \caption{Each figure shows a $z=0$ slice capturing the flow box setup of one of Corollaries~\ref{extend3}--\ref{extend8}.  Plus signs are positioned on the side from which the contact structure dominates the horizontal foliation.  Dashes, for instance along $U$, show where the confoliation slope is 0.  The smooth foliation acts as a transport mechanism for contact structures in the direction shown by the arrows.}
   \label{extend}
\end{figure}

\begin{cor}[Figure~\ref{extend}a]\label{extend3} Let $U$ be a regular neighborhood in $F$ of the union of $x=\pm 1$ and $z=\pm 1$.  Let $\eta$ be a $C^k$-confoliation  with $k \ge 1$ defined in a neighborhood $V$ of $y=1$ in $F$ by a 1-form $$dz - a(x,y,z) dx$$ with
\begin{enumerate}
\item  $a(x,y,z)=0$ on $V\cap \overline{U}$, and 
\item  $a(x,y,z)>0$ and $\dfrac {\partial a}{\partial y}(x,y,z) > 0$ on $V\setminus \overline{U}$.
\end{enumerate}
Then $\eta$ can be extended to a $C^k$-confoliation $\hat\eta$ on $F$ that agrees with the horizontal foliation of $F$ in $\overline{U}$, is contact on the complement of $\overline{U}$, and strictly dominates, when viewed from outside $F$, the horizontal foliation on the complement in $y=0$ of $\overline{U}$. 
\end{cor}

\begin{proof}
It is enough to extend $a(x,y,z)$ along each interval $\{x_0\} \times [0,1] \times \{z_0\}$ to a $C^k$ function $\hat{a}(x,y,z)$ such that 

\begin{enumerate}
\item  $\hat{a}(x,y,z)=0$ in $\overline{U}$, and 
\item  $\hat{a}(x,y,z)>0$ and $\dfrac {\partial \hat{a}}{\partial y}(x,y,z) > 0$ outside $\overline{U}$.
\end{enumerate} 
\end{proof}

\begin{cor}[Figure~\ref{extend}b]\label{extend4} Let $U$ be a  regular neighborhood in $F$ of the union of $x=\pm 1$ and $z=\pm 1$.  Let $\eta$ be a $C^k$-confoliation  with $k \ge 1$ defined in a neighborhood $V$ of $y=1$ in $F$ by a 1-form $$dz - a(x,y,z) dx$$ with
\begin{enumerate}
\item  $a(x,y,z)=0$ on $V\cap \overline{U}$, and 
\item  $a(x,y,z)>0$ and $\dfrac {\partial a}{\partial y}(x,y,z) > 0$ on $V\setminus \overline{U}$.
\end{enumerate}
Then $\eta$ can be extended to a $C^k$-confoliation $\hat\eta$ on $F$ that agrees with the horizontal foliation of $F$ in $\overline{U}$, is contact on the complement of $\overline{U}$ in the interior of $F$, and  smoothly agrees with the horizontal foliation at $y=0$. 
\end{cor}

\begin{proof} Proceed as in the proof of  Corollary~\ref{extend3} except insist that $$\hat{a}(x,0,z)\equiv 0.$$
\end{proof}

\begin{cor}[Figure~\ref{extend}c]\label{extend5} Let $U$ be a  regular neighborhood in $F$ of the union of $x=\pm 1$ and $z=\pm 1$.  Let $\eta$ be a $C^k$-confoliation with $k \ge 1$ defined in a neighborhood $V$ of the union of $y=0$ and $y=1$ in $F$. Suppose that when viewed from inside $F$, $\eta$ dominates the horizontal foliation along the vertical faces given by
$y=0$ and $y=1$, with strict domination  in the complement of $\overline{U}$. Then $\eta$ can be extended to a $C^k$-confoliation $\hat\eta$ on $F$ that agrees with the horizontal foliation of $F$ in $\overline{U}$, and  is contact on the complement of $\overline{U}$.
\end{cor}

\begin{proof}
Decompose $F$ as a union of two flow boxes diffeomorphic to $F$ by cutting open along the plane $y=1/2$. Apply Corollary~\ref{extend4} to each of the resulting flow boxes.
\end{proof}

The point of ``smoothly agrees with''  in the next corollary is that flow boxes are brick-like objects that, when sensibly glued together, should define a smooth confoliation.  Thus we require smooth convergence of the confoliation to horizontal at $y=0$.

\begin{cor}[Figure~\ref{extend}d]\label{extend6} Let $U$ be a  regular neighborhood in $F$ of the union of $x=\pm 1$ and $z=\pm 1$.  Let $\eta$ be a $C^k$-confoliation  with $k \ge 1$ defined in a neighborhood $V$ of $y=1$ in $F$ by a 1-form $$dz - a(x,y,z) dx$$ with
\begin{enumerate}
\item  $a(x,y,z)=0$ on $V\cap \overline{U}$, and 
\item  $a(x,y,z)>0$ and $\dfrac {\partial a}{\partial y}(x,y,z) > 0$ on $V\setminus \overline{U}$.
\end{enumerate}
Then $\eta$ can be extended to a $C^k$-confoliation $\hat\eta$ on $H$ that agrees with the horizontal foliation of $H$ in $\overline{U}$, is contact on the complement of $\overline{U}$ in the interior of $H$, smoothly agrees with the horizontal foliation at $y=0$, and dominates the horizontal foliation in the complement of $\overline{U}$  along the lines 
$$y = x  -1/2 \mbox{      and       } y = -1/2-x. $$  
\end{cor}

\begin{proof} This follows from Corollary~\ref{extend4}. 
\end{proof}

\begin{cor}[Figure~\ref{extend}e]\label{extend7}
Let $U$ be a  regular neighborhood in $G$   of the union of $z=\pm 1$ and $\Delta^{(0)}\times [-1,1]$.  
Let $S$ denote the complement in $\partial_v G$ of the 2-cell given by $y=1$. Let $\eta$ be a $C^k$-confoliation with $k \ge 1$ defined in a neighborhood $V$ of $S$ in $G$ such that 
\begin{enumerate}
\item $\eta$  is contact on $V\setminus \overline{U}$, 
\item $\eta$ agrees with the horizontal foliation on $V\cap\overline{U}$, and 
\item $\eta$ strictly dominates the horizontal foliation along $S\setminus \overline{U}$, when viewed from inside $G$.
\end{enumerate}
Then $\eta$ can be extended to a $C^k$-confoliation $\hat\eta$ on $G$ that agrees with the horizontal foliation of $G$ in $\overline{U}$, is contact on the complement of $\overline{U}$ in the interior of $G$, and, when viewed from outside $G$, dominates the horizontal foliation in the complement of $\overline{U}$  along the line $y=1$.
\end{cor}

\begin{proof}
Let $\alpha$ be a smooth arc  properly embedded in $\Delta$ which connects the vertices $(3/2,1)$ and $(0,-1/2)$ and is not tangent to a side of the triangle $\Delta$ at its endpoints. Let $R=\alpha\times [-1,1]$, a vertical rectangle in $G$. Decompose $G$ along $R$ into two flow boxes as $G=G'\cup F'$, where $G'$ is diffeomorphic to $G$ and $F'$ is diffeomorphic to $F$. First apply Corollary~\ref{extend6} to $G'$ and then apply Corollary~\ref{extend3} to $F'$.
\end{proof}

\begin{cor}[Figure~\ref{extend}f]\label{extend8}
Let $U$ be a  regular neighborhood in $G$   of the union of $z=\pm 1$ and $\Delta^{(0)}\times [-1,1]$.  Let $\eta$ be a $C^k$-confoliation with $k \ge 1$ defined in a neighborhood $V$ of $\partial_v G$ in $G$ such that 
\begin{enumerate}
\item $\eta$  is contact on $V\setminus \overline{U}$,
\item $\eta$ agrees with the horizontal foliation on $V\cap\overline{U}$, and 
\item $\eta$ strictly dominates the horizontal foliation along $\partial_v G \setminus \overline{U}$, when viewed from inside $G$.
 \end{enumerate}
 Then $\eta$ can be extended to a $C^k$-confoliation $\hat\eta$ on $G$ that agrees with the horizontal foliation of $G$ in $\overline{U}$ and is contact on the complement of $\overline{U}$ in the interior of $G$.

\end{cor}

\begin{proof}
This time we cut  $G$ open along two vertical rectangles and consider instead the resulting union of flow boxes. This time we cut $G$ open along two disjoint vertical rectangles to obtain
$$G=G'\cup F'\cup F^{''},$$
where $G'$ is diffeomorphic to $G$ and each of $F'$ and $F^{''}$ is diffeomorphic to $F$. Then apply Corollary~\ref{extend6} to $G'$. Finally, apply Corollary~\ref{extend3} to each of $F'$ and $F^{''}$.
\end{proof}

We now consider a case where the initial confoliation is defined on the entire vertical boundary of a solid cylinder.
Let $\mathcal L$ be a smooth 1-dimensional foliation on the cylinder $S^1 \times I$ such that the boundary components of $S^1\times I$ are leaves of $\mathcal L$.  Let $h:I\to I$ be the holonomy map of $\mathcal L$.

\begin{prop}\label{cylinderextend} If $h'(z)<0$ for $z\in(0,1)$, there is a confoliation on $D^2 \times I$ that is contact on $D^2 \times (0,1)$, tangent  on $D^2 \times \partial I$, has characteristic foliation $\mathcal L$, and is everywhere transverse to $\partial/\partial z$.
\end{prop}

\begin{proof} This will follow from Lemmas~\ref{extension1}--\ref{extension3}.
\end{proof}

\begin{lemma}\label{extension1} If $h'(z)<0$ for $z\in(0,1)$, then there exists a smooth 1-dimensional foliation $\mathcal K$ of $S^1 \times I$ with the same holonomy map $h:I\to I$ such that for every $(\theta, z)\in S^1 \times (0,1)$, the slope, $s(\theta,z)$, of $\mathcal K$ at $(\theta,z)$ is negative.
\end{lemma}

\begin{proof} As a first approximation, let $\mathcal K_1$ be the foliation of $[0,2\pi] \times I$ given by connecting each point $(0,z)$ to $(2\pi,h(z))$ by a straight line.  Next create $S^1 \times I$ by identifying $(0,z)$ and $(2\pi,z)$ for $z\in I$.  Let $\mathcal K_2$ denote the image of $\mathcal K_1$ in $S^1\times I$.  This has the desired properties, except that $\mathcal K_2$ is not smooth along $\{0\} \times I$. Carefully rounding these corners (see for example, Lemma 4.7 of \cite{Milnor}), yields the desired smooth foliation.
\end{proof}

\begin{lemma} \label{extension2} There is a diffeomorphism, $F$, of $D^2 \times I$ that is the identity map on $D^2 \times \partial I$ and takes $\mathcal L$ to $\mathcal K$.
\end{lemma}

\begin{proof} Let $\theta = 0$ be a base point for $S^1$ so that the holonomy map for each foliation is $h:\{0\} \times I \to \{0\} \times I$.  Let $f$ be the diffeomorphism of $S^1 \times I$ such that $f$ restricts to the identity map on $\{0\} \times I$, preserves the $S^1$ coordinate, and maps $\mathcal L$ to $\mathcal K$.

Define $F:D^2 \times I\to D^2 \times I$ by $F(r,\theta,z)=(r, t(r)f(\theta,z) + (1-t(r))(\theta,z))$ where $t$ is a diffeomorphism of the interval smoothly damped at the endpoints. 
\end{proof}

\begin{lemma} \label{extension3} There is a confoliation $\xi$ on $D^2 \times I$ that is contact on $D^2 \times (0,1)$, tangent on $D^2 \times \partial I$, has characteristic foliation $\mathcal K$, and is everywhere transverse to $\partial/\partial z$.
\end{lemma}

\begin{proof} Using cylindrical coordinates, define $\alpha = dz - r^2s(\theta,z)d\theta$.  Then $d\alpha=-2rs(\theta,z)drd\theta-r^2s_z(\theta,z)dzd\theta$, from which it follows that 
$$\alpha \wedge d\alpha=-2s(\theta,z)rdrd\theta dz.$$
Setting $\xi=\text{ker}(\alpha)$ has the desired properties.
\end{proof}

\section{The inductive construction of $\xi$.}\label{theconstructionsection}

In this section, we show how a $C^{1,0}$-foliation can be used to propagate a contact structure across $M$.  Before describing this procedure, we highlight the role of  smoothness in the approach used by Eliashberg and Thurston in performing this propagation.

First,  recall Lemma~\ref{ETlemma}.  Roughly speaking,  a clever choice of foliation coordinates permits a confoliation along a Legendrian curve to be described by a monotone function.  By taking advantage of a beginning contact zone,  such a function can be approximated by a strictly monotone function, thereby creating a larger contact zone.  This argument can be repeated on overlapping regions covering the manifold.  The issue is whether strict monotonicity attained on a given region can be preserved under subsequent approximations.  This is precisely where smoothness of the foliation becomes important,  guaranteeing that derivatives are globally defined and continuous, and thus allowing one to preserve strict monotonicity under subsequent approximations. 

We circumvent the issue of monotonicity with carefully chosen, minimally overlapping, flow boxes and a more discrete propagation technique, which we now describe.

Recall Definition~\ref{dominate} and the slope convention chosen in Definition~\ref{flowboxdefn}. Let $F$ be a flow box and let $R$ be a rectangle in $\partial_v F$. Let  $\chi_{\xi_1}$ and $ \chi_{\xi_2}$ denote the characteristic foliations induced on $R$  by two 2-plane fields $\xi_1$ and $\xi_2$ defined in a neighborhood of $R$ and positively transverse to $\Phi$. We  write $\chi_{\xi_1}<_p\chi_{\xi_2}$ if the unit vector tangent to $\chi_{\xi_1}$ at $p$ has slope less than the unit vector tangent to $\chi_{\xi_2}$ at $p$.  Similarly, we write $\chi_{\xi_1}=_p\chi_{\xi_2}$ if the unit vector tangent to $\chi_{\xi_1}$ at $p$ has slope equal to  the unit vector tangent to $\chi_{\xi_2}$ at $p$. We write $\chi_{\xi_1}<\chi_{\xi_2}$ if $\chi_{\xi_1}<_p\chi_{\xi_2}$ for all $p\in \text{int}(R)$. Similarly, we write $\chi_{\xi_1}\le\chi_{\xi_2}$ if $\chi_{\xi_1}\le_p\chi_{\xi_2}$ for all $p\in \text{int}(R)$.
\begin{definition}

Let $M$ be a closed oriented 3-manifold with smooth flow $\Phi$. Suppose that $M$ decomposes as  a union 
$$M=V\cup W,$$
where $V$ and $W$ are smooth 3-manifolds, possibly with corners, and $\partial V=\partial W$. We say that this decomposition  is {\sl compatible with the flow $\Phi$} if $\partial V$ (and hence $\partial W$) decomposes as a union of compact subsurfaces $\partial_v V\cup \partial_h V$, where $\partial_v V$ is a union of flow segments of $\Phi$ and, $\partial_h V$ is transverse to $\Phi$.   In the presence of a foliation transverse to the flow, the notation $\partial_h V$ will be used for the portion of $\partial V$ tangent to the foliation.
\end{definition}

\begin{definition}\label{compatible1} Let $M$ be a closed oriented 3-manifold with smooth flow $\Phi$. Suppose that $M$ can be expressed as  a union 
$$M=V\cup W,$$
where $V$ and $W$ are smooth 3-manifolds, possibly with corners, such that  $\partial V=\partial W$.  Suppose also that
this decomposition is compatible with $\Phi$, that $V$ admits a smooth contact structure  $\xi_V$, and that $W$ admits a $C^{1,0}$-foliation $\mathcal F_W$.  We say that $(V,\xi_V)$ is {\sl $\Phi$-compatible} with $(W,\mathcal F_W)$, and that $M$ admits a positive $(\xi_V,\mathcal F_W,\Phi)$ decomposition, if the following are satisfied:
\begin{enumerate}
\item $\xi_V$ and $\mathcal F_W$ are (positively) transverse to $\Phi$ on their domains of definition,
\item $\xi_V$ is tangent to $\partial_h V$, and
\item $\chi_{\xi_V}<\chi_{T\mathcal F_W}$ on the interior of $\partial_v V$, when viewed from outside $V$.
\end{enumerate}

\end{definition}

The main result of the section is Theorem~\ref{main1}.  The starting point is a transitive flow box decomposition $M= V\cup F_1 \cup \dots \cup F_n$ for $\mathcal F$.  Set $V_0 = V$ and $V_i = V\cup F_1\cup...\cup F_i$ for each $i, 1\le i\le n$. For each $i, 0\le i\le n$, set  $W_i = M\setminus \text{int}(V_i)$. Thus,  for $0\le i \le n$, $\mathcal F$ is compatible with $(V_i,W_i,\Phi)$.  When $i=0$, $\partial_h V_i = \emptyset$, and when $i=n$, $V_n=M$ and hence $\partial V_n=\emptyset$. 

\begin{lemma} \label{smoooth} After possibly blowing up finitely many leaves of $\mathcal F$, we may  assume that
\begin{enumerate}
\item $\mathcal F$ is smooth in a neighborhood of the horizontal faces $\partial_h F_i$, and
\item $\mathcal F$ is smooth in a neighborhood of the vertical edges of each $F_i$.
\end{enumerate}
\end{lemma}

\begin{proof} 
We will define $\xi$ inductively and in a piecewise fashion over each  $F_i$. To guarantee that the resulting confoliation is everywhere smooth,   we add a {\sl smoothly foliated collar} about the horizontal boundary and the vertical 1-simplices of each $F_i$ as follows.
Let $L_1, L_2, \dots, L_m$ be the collection of leaves of $\mathcal F$ that contain the horizontal boundaries of all $F_i$.  Modify the original foliation $\mathcal F$ by thickening each of the leaves $L_j$.  Thus each $L_j$ is replaced with a product of leaves $L_j \times [-1,1]$.  The thickening should be performed so that for each $i$, if a component of $\partial_h F_i$ was originally contained in $L_j$, then it is now contained in $L_j \times \{0\}$. As noted in Lemma~\ref{blowup}, we may assume the interior of each thickening, $L_j\times (-1,1)$, is smoothly immersed. Note that $M= V\cup F_1 \cup \dots \cup F_n$ remains a transitive flow box decomposition of $M$ with respect to this new foliation. Moreover, $\mathcal F \cap F_i$ is smooth in a neighourhood of $\partial_h F_i$ for each $i$.

The vertical edges of the $F_i$ are smooth transverse arcs, thus $\mathcal F$ can be smoothed on $D^2 \times I$ product neighborhoods of these edges.
\end{proof}

 Now we fix a preferred regular neighborhood of the union of the horizontal 2-cells and the vertical 1-cells of the flow box decomposition. We do this as follows. For each $i, 1\le i \le n$, choose a regular neighborhood of $ \partial_h F_i$ which is contained in the thickening $L_i\times (-1,1)$, and  choose a regular neighborhood of the vertical 1-cells of $F_i$on which $\mathcal F$ is smooth. Choose these neighborhoods  so that the union $U$ of these neighborhoods is a regular neighborhood of the union of all horizontal 2-cells and all vertical 1-cells of the flow box decomposition. We refer to both this preferred neighborhood $U$ and the restriction of $\mathcal F$ to $U$ as the {\sl smoothly foliated collar}.

Next, we modify Definition~\ref{compatible1} slightly to account for the smoothly foliated collar.

\begin{definition}\label{compatible} Let $0\le i<n$ and let $\xi_i$ be a smooth confoliation on  $V_i$ such that $\xi_i$ is tangent to $\partial_h V_i$ and everywhere transverse to $\partial_v V_i$. (Note that by smoothness, we may consider $\xi_i$ to be defined on an open neighborhood of $ V_i$.) We say the smooth confoliation $\xi_i$ {\sl dominates} $\mathcal F$ if $\chi_{\xi_i} \ge \chi_{\mathcal F}$ on $\partial_v V_i$ when viewed from outside $V_i$, with equality permitted only on the closure, $\overline{U}$, of the smoothly foliated collar. Use {\sl strictly dominates} if the inequality is strict.  

Let $\xi_i$ be a smooth confoliation defined on $V_i$.  We say that $(V_i,\xi_i)$ is {\sl  smoothly $\Phi$-compatible} with $(W_i,\mathcal F)$ if the following are satisfied:
\begin{enumerate}
\item $\xi_i$ and $\mathcal F_W$ are (positively) transverse to $\Phi$ on their domains of definition,
\item $\xi_i$ is tangent to $\partial_h V_i$ 
\item $\xi_i = T\mathcal F$ on $\overline{U}\cap V_i$,
\item $\xi_i$ is a contact structure on  $V_i\setminus U$, and
\item $\xi_i$ dominates $\mathcal F$ (on $\partial_v V_i$).
\end{enumerate}
\end{definition}

\begin{prop} \label{inductthis}Let $M = V \cup F_1 \cup \dots \cup F_n$ be a transitive flow box decomposition of $M$, and let $\xi_i$, $i\ge 0$, be a smooth confoliation defined on $V_i$ that is smoothly $\Phi$-compatible with $(W_i, \mathcal F)$.  Then there is a smooth confoliation $\xi_{i+1}$ defined on $V_{i+1}$ that is smoothly $\Phi$-compatible with $(W_{i+1}, \mathcal F)$ and restricts to $\xi_i$ on $V_i$.
\end{prop}

\begin{proof} Since any polygon admits a triangulation, any transitive flow box decomposition can be chosen to consist only of flow boxes  diffeomorphic to the flow box $G=\Delta\times [0,1]$ defined in Section~\ref{basicsection}. In particular, we may assume that $F_{i+1}$ is diffeomorphic to $G$. Such a diffeomorphism preserves the foliation and flow directions, that is, slopes $0$ and $\infty$;  thus we will make slope comparisons and approximations without reference to the change of coordinates.

By hypothesis, $\xi$ strictly dominates $\mathcal F$ on $X=\partial_v V_i \setminus U$.  By compactness, there exists $\epsilon>0$ such that 
$$\text{slope}(\xi_i)-\text{slope}(\mathcal F)> 3\epsilon$$ on $X\cap \partial_v F_{i+1}$.

By Proposition~\ref{smoothapprox}, we may approximate the restriction of $\mathcal F$ to $F_{i+1}$ by a smooth foliation $\tilde{\mathcal  F}$ such that $(\tilde{\mathcal F}, T\tilde{\mathcal F})$ is $C^0$-close to $(\mathcal F,T\mathcal F)$.  Choose this approximation so that $\xi_i$ dominates $\tilde{\mathcal F}$ on $\partial_v(F_{i+1}) \cap V_i$, and so that 
$$|\text{slope}(\tilde{\mathcal F})-\text{slope}(\mathcal F)|<\epsilon$$ on all of $\partial_v F_{i+1}$.  It follows that 
$$\text{slope}(\xi_i)-\text{slope}(\tilde{\mathcal F}) > 2\epsilon$$ on $X\cap \partial_v F_{i+1}$. Choose smooth coordinates $(x,y,z)$ on $F_{i+1}$ so that the leaves of $\tilde{\mathcal F}$ are horizontal, given by $z=$ constant. (Note that although this change of coordinates might change slope values, it doesn't affect the relative values of slopes.)

Next consider the number of  2-cells contained in  $\partial_v(F_{i+1}) \cap V_i$.  From the definition of a transitive flow box decomposition there is at least one.  Depending on whether there are exactly one, two, or three such 2-cells, apply the corresponding Corollary~\ref{extend6}, \ref{extend7}, or \ref{extend8} to smoothly extend $\xi_i$ across $F_{i+1}$ and call the resultant confoliation $\xi_{i+1}$.  Smoothness of the glued confoliation is assured by the construction of the confoliations near the boundaries of their flow boxes.

In constructing the extensions of Corollaries~\ref{extend5}--\ref{extend8} the starting point is a 1-form on $\partial_v F_{i+1}$ given by  $dz-a(x,y,z)dx$ with $a(x,y,z)>0$ on $ X \cap \partial_v F_{i+1}$.  This is then extended across $F_{i+1}$ while keeping $a(x,y,z)>0$ and also $\frac{\partial a}{\partial y}(x,y,z)>0$.  Given any $\delta>0$, these constructions can be performed while keeping the change in $a(x,y,z)$ along Legendrian curves less than $\delta$.  In other words, the extension $\xi_{i+1}$ can be chosen so that the change in $\text{slope}(\xi_{i+1})$ along Legendrian curves is less than $\epsilon$.  Thus 
$$\text{slope}(\xi_{i+1})-\text{slope}(\tilde{\mathcal F}) > \epsilon$$ on $\partial_v V_{i+1} -(N_h\cup N_v)$, and consequently, on this set we also have
$$\text{slope}(\xi_{i+1})-\text{slope}({\mathcal F}) > 0.$$
\end{proof}

\begin{cor}\label{inductthiscor}
Let $M = V \cup F_1 \cup \dots \cup F_n$ be a transitive flow box decomposition of $M$, and let $\xi_i$, $i\ge 0$, be a smooth confoliation defined on $V_i$ that is smoothly $\Phi$-compatible with $(W_i, \mathcal F)$ and lies within $\epsilon$ of $\mathcal F$ on the intersection of the domain of  $\mathcal F$ with  $V_i$.  Then there is a smooth confoliation $\xi_{i+1}$ defined on $V_{i+1}$ that is  smoothly $\Phi$-compatible with $(W_{i+1}, \mathcal F)$, restricts to $\xi_i$ on $V_i$, and lies within $3\epsilon$ of $T\mathcal F$ on the intersection of the domain of  $\mathcal F$ with  $V_{i+1}$.
\end{cor}

\begin{proof}
It follows from Proposition~\ref{smoothapprox} that a smooth foliation $\tilde{\mathcal F}$ may be chosen on $F_{i+1}$  so that $T\tilde{\mathcal F}$ is within $\epsilon$ of $T\mathcal F$.  Restricting attention to $\partial_v V\cap F_{i+1}$, $T\tilde{\mathcal F}$ lies within $2\epsilon$ of $\xi_i$.  From the proof of Proposition~\ref{inductthis}, $\xi_{i+1}$ is constructed to be as close or closer  to $T\tilde{\mathcal F}$ on $F_{i+1}$ than $\xi_i$ is on $\partial_v V\cap F_{i+1}$.  Thus $\xi_{i+1}$ is within $2\epsilon$ of $T\tilde{\mathcal F}$, and hence within $3\epsilon$ of $T\mathcal F$.
\end{proof}

\begin{cor}\label{almost} Let $M = V \cup F_1 \cup \dots \cup F_n$ be a transitive flow box decomposition of $M$, and let $\xi_0$  be a smooth contact structure defined on $V$ that is compatible with $(W, \mathcal F)$.  Then there is a smooth contact structure  defined on $M$ that agrees with $\xi_0$ in $V$ and is $\Phi$-close to $\mathcal F$ on  $W$.
\end{cor}

\begin{proof} Inductively applying Proposition~\ref{inductthis} produces a transitive smooth confoliation that is $\Phi$-close  to $\mathcal F$ on $W$ and agrees with $\xi_0$ in $V$.  By Proposition~2.8.1 of \cite{ET} (see also \cite{Et2}), this transitive smooth confoliation can be smoothly deformed into a smooth contact structure. Thus there is a smooth contact structure  defined on $M$ that agrees with $\xi_0$ in $V$ and is $\Phi$-close to $\mathcal F$ on  $W$.
\end{proof}

\begin{cor}
Let $M = V \cup F_1 \cup \dots \cup F_n$ be a transitive flow box decomposition of $M$, and let $\xi_0$  be a smooth confoliation defined on $V$ that is $\Phi$-compatible with $(W, \mathcal F)$. Then there is a smooth contact structure  defined on $M$ that agrees with $\xi_0$ in $V$ and is $\Phi$-close to $\mathcal F$ on  $W$. If $\xi_0$ is within $\epsilon$ of $T\mathcal F$ along $\partial V=\partial W$, then $\xi$ can be chosen to lie within $f(n)\epsilon$ of $T\mathcal F$ in $W$, for some positive function $f$.
\end{cor}

\begin{proof}
This follows immediately from Corollary~\ref{inductthiscor}.
\end{proof}

\begin{lemma}\label{phitoC0} Let $M$ be a closed oriented 3-manifold. Suppose that $M$ can be expressed as  a union 
$$M=V\cup W,$$
such that $\partial V=\partial W$, $V$ admits a smooth contact structure $\xi_0$, $W$ admits a $C^{1,0}$-foliation $\mathcal F$, and $(V,\xi_0)$ and $(W,\mathcal F)$ are $\Phi$-compatible.  Then $\xi_0$ can be modified in an arbitrarily small collar neighborhood of $\partial V$ so that the restriction of $\xi_0$ to $\partial_v V$ lies arbitrarily $C^0$-close to $T\mathcal F$.
\end{lemma}

\begin{proof} Fix any $\epsilon>0$. Let  $X=(\partial_v V) - (N_h\cup N_v).$ Let $X\times (-\delta,0]$ be a collar neighborhood of $X$ in $V$, with $X=X\times \{0\}$. Pick any smooth line field $l$ on $X$ which is tangent to $\mathcal F$ along the smoothly foliated collar, is dominated by the projection of  $\xi_0|_{X\times \{ -\delta\}}$ to $X$, dominates the projection of $\mathcal F|_{X\times \{0\}}$ to $X$, and lies with $\epsilon$ of $\mathcal F|_{X\times \{0\}}$. Replace the restriction of $\xi_0$ to $X\times  [-\delta,0]$ by a  straight line homotopy between $\xi_0$ restricted to $X\times \{-\delta\}$ and $l$, damped to fit smoothly with the original $\xi_0$ as defined on the complement of the collar.
\end{proof}

\smallskip
\begin{thm}~\label{main1} If $M$ admits a positive $(\xi_V,\mathcal F_W,\Phi)$ decomposition, then $M$ admits a smooth positive contact structure   $\xi_+$  which agrees with $\xi_V$ on $V$ and is  $\Phi$-close to $\mathcal F_W$ on $W$.  Moreover, $M$ also admits a smooth contact structure which agrees with $\xi_0$ on the complement in $V$ of a collar neighborhood of $\partial V$ and is  arbitrarily  $C^0$-close to $\mathcal F$ on $W$. The analogous result holds if $M$ admits a negative $(\xi_{V'},\mathcal F_{W'},\Phi)$ decomposition. If $M$ admits both a positive $(\xi_V,\mathcal F_W,\Phi)$ decomposition and a negative $(\xi_{V'},\mathcal F_{W'},\Phi)$ decomposition, then these contact structures   $(M,\xi_+)$ and $(-M,\xi_-)$ are weakly symplectically fillable and universally tight.  
\end{thm}

\begin{proof} This follows immediately from Proposition~\ref{transitiveflowbox}, Corollary~\ref{almost}, and Lemma~\ref{phitoC0}.
\end{proof}

\section{Attracting holonomy}\label{attracting holonomy}

This section gives a generalization of the Eliashberg-Thurston result on perturbing foliations in a neighborhood of a curve in a leaf with sometimes attracting holonomy (Proposition 2.5.1, \cite{ET}) to the larger class of $C^{1,0}$-foliations  with holonomy containing a {\sl contracting interval}.  Holonomy with a contracting interval will be defined below in terms of the existence of a particular smooth submanifold with corners.  Holonomy with a contracting interval is a weaker condition than holonomy with an attracting leaf.

Let $P$ be the prism in $\mathbb R^3$ given by $|x|\le 1$, $|y|\le 1$ and $|z|\le y/2 +3/2$.  The slanted top and bottom of $P$ is denoted by $\partial_h P$, and the rest of $\partial P$ is denoted by $\partial_v P$.  Let $V$ be the solid torus given by identifying each pair of points $(x,-1,z)$ and $(x,1,z)$ where both are in $P$.  Then $\partial_h V$ is defined to be the image of $\partial_h P$, and $\partial_v V$ is defined to be the portion of $\partial V$ that is in the image of $\partial_v P$.

\begin{definition}\label{contracting}

The holonomy of a $C^{1,0}$-foliation $\mathcal F$ transverse to a flow $\Phi$ in a 3-manifold $M$ has a {\sl contracting interval} if there is a subset of $M$ diffeomorphic to $V$, called an {\sl attracting neighborhood}, such that

\begin{enumerate}

\item $\partial_h V$ is mapped into leaves of $\mathcal F$ and

\item vertical intervals in $V$ are mapped to flow lines of $\Phi$.
\end{enumerate}
\end{definition}

\begin{thm}\label{attracting} Let $\mathcal F$ be a taut oriented $C^{1,0}$-foliation transverse to a flow $\Phi$.  If $V$ is a disjoint union of  attracting neighborhoods, and $\mathcal F$ is $V$-transitive, then $\mathcal F$ can be $\Phi$-approximated both by a smooth positive contact structure   $\xi_+$ and by a smooth negative contact structure   $\xi_-$.     These contact structures $(M,\xi_+)$ and $(-M,\xi_-)$ are weakly symplectically fillable and universally tight. 
\end{thm}

\begin{proof}  In light of Theorem~\ref{main1}, it is enough to construct a smooth confoliation $\xi$ on $V$ such that $(V, \xi)$ is compatible with $(W,\mathcal F)$, in the sense of Definition~\ref{compatible}. It is sufficient to consider the case that $\xi$ is positive. The construction of the desired $\xi$ is the same on each component of $V$, so we will treat $V$ as if it is connected. Let $W$ denote the closure of the complement of $V$ in $M$.

 \begin{figure}[htbp] %  figure placement: here, top, bottom, or page
   \centering
   \includegraphics[width=5in]{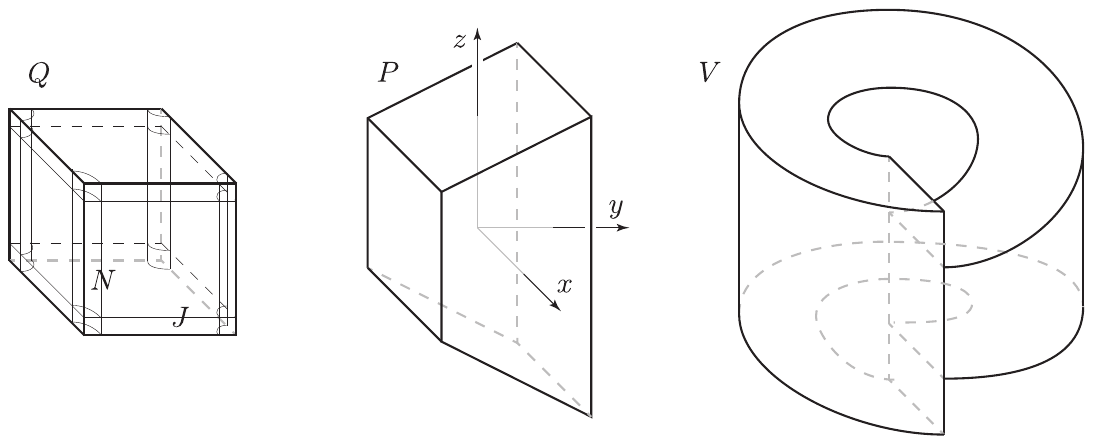} 
   \caption{Foliated neighborhoods of horizontal faces and vertical edges are shown only in $Q$ though they exist in $P$ and $V=P/\sim$ as well.}
   \label{QPV}
\end{figure}

Consider the transformation of the prism $P$ to the cube $$Q=[-1,1] \times [-1,1] \times  [-1,1]$$ that fixes $(x,y)$ and linearly scales the $z$-coordinate. View $V$ as the quotient $$V=Q/{\sim}$$ of $Q$ obtained by identifying
$$(x,-1,z)\sim (x,1,z/2).$$
Notice that $\partial_h Q$ is mapped into leaves of $\mathcal F$, and vertical intervals in $Q$ are mapped to flow lines of $\Phi$. Moreover, we may assume that the original parametrization of $P$ was chosen so that $\mathcal F$ meets the $y=\pm 1$ sides of $\partial_v Q$ in horizontal lines. Since $\mathcal F$ is $C^0$, the leaves of the foliation meet each of the  $x=\pm 1$ sides in a continuous family of smooth graphs. 

To facilitate smooth gluings, thicken the leaf or leaves of $\mathcal F$ which meet $\partial_h Q$. Fix $0< \epsilon< 1/4$.  Choose the thickening of the leaves of $\mathcal F$  intersecting  $\partial_h V$ to replace the leaves $z=\pm 1$ in $Q$ with a disjoint union, $J$, of $I$-bundles in $Q$, with $J$ containing $[-1,1]\times[-1,1]\times [-1,-1+\epsilon]$ and $[-1,1]\times[-1,1]\times [1-\epsilon,1]$ as components. Assume also that $\mathcal F$ meets an $\epsilon$-neighborhood, $N$, of the quotient of the vertical edges of $Q$ in a smooth horizontal foliation. Let $U_Q = \text{int} (J\cup N)$;  $U_Q$ is  an $\epsilon$-neighbourhood of the union of $\partial_h Q$ with the 1-cells of $\partial_v Q$. We assume also  that the pull-back of $\mathcal F$  to $\overline{U_Q}$ is horizontal. We will abuse notation and let $U_Q$ also refer to the projection of $U_Q$ to $M$.

Since $Q$ is a flow box it is amenable to the constructions of Section~\ref{basicsection}.  We will $\Phi$-approximate $\mathcal F$ in $Q$, and thus in $V$, by a $C^{\infty}$-confoliation $\xi_0$ which smoothly respects both the identification $\sim$ and the gluing of $(V,\xi_0)$ and $(W,\mathcal F)$  along $\partial_h V=\partial_h W$. 

The  confoliation $\xi_0$ will be chosen to agree with $\mathcal F$ on $U_Q$ and to be a contact structure on $Q\setminus \overline{U_Q}$ which, when viewed from outside $Q$, is strictly dominated by $\mathcal F$ on the $y=-1$ side of $\partial_v Q$ and strictly dominates $\mathcal F$ on the remaining three sides of $\partial_v Q$.  See Figure~\ref{thecube}. Since the interior of the $y=-1$ side lies in the interior of $V$, the resulting $(V,\xi_0)$ will be compatible with $(W,\mathcal F)$ in the sense of Definition~\ref{compatible}. 

 \begin{figure}[htbp] %  figure placement: here, top, bottom, or page
   \centering
   \includegraphics[width=5in]{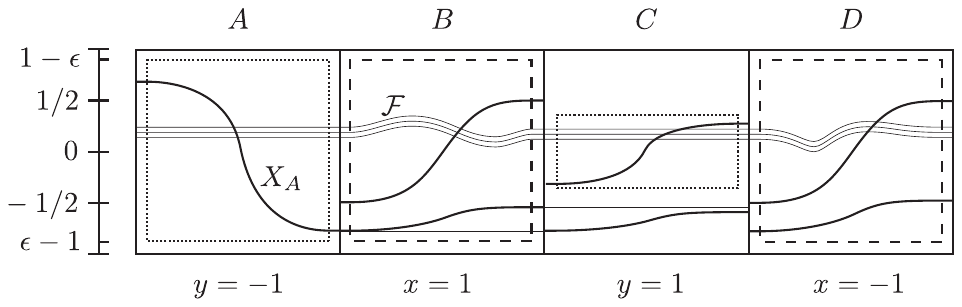} 
   \caption{Some of the choices of flow lines of the vector fields $X_A,...$ and their relationship to $\mathcal F$ are shown.  Not enough detail is drawn to show that the holonomy given by flowing from left to right, that is counterclockwise about $\partial_v Q$, is decreasing.}
   \label{thecube}
\end{figure}

As a first step in constructing $\xi_0$, we define $\xi_0$ along $\partial_v Q$.  The vertical boundary $\partial_v Q$ consists of four vertical faces. Let $A$ denote the face $y=-1$, let $B$ denote the face $x=1$, let $C$ denote the face $y=1$, and let $D$ denote the face $x=-1$. 

We construct $\xi_0$ by first specifying  smooth unit vector fields $X_A, X_B, X_C,$ and $X_D$ along the faces $A, B, C, $ and $D$ respectively, and then declaring $\xi_0$ along $\partial_v Q$  to be the 2-plane field which is normal to $\partial_v Q$ and contains the corresponding tangent vector $X_A, X_B, X_C,$ or $ X_D$. We will choose $X_A, X_B, X_C,$ and $X_D$ to be horizontal on $\overline{U_Q}$ and hence in a neighborhood of the vertical 1-simplices of $\partial_v Q$; in particular, the 2-plane field $\xi_0$ will therefore be well-defined   on the vertical edges.

Begin by defining the vector field $X_B$.  Choose $X_B=X_B(1,y,z)$  to be a smooth unit vector field which satisfies the following
\begin{enumerate}
\item $X_B$ dominates $\mathcal F$,
\item $X_B$ has positive slope when both $y$ and $z$ lie in $(-1+\epsilon, 1-\epsilon)$, and
\item $X_B =\partial/\partial y$ when   $y$ or $z$ lies in $ [-1,-1+\epsilon]\cup [1-\epsilon,1]$.
\end{enumerate}

Let $\Psi_B$ denote the flow generated by $X_B$.  Abusing notation a bit, we denote by $\Psi_B(y,z)$ the intersection of of the flow line of $\Psi_B$ that starts at $(1,-1,z)$ with $\{(1,y)\}\times[-1,1]$. (We  use this notation  when referring to all flows in this section.)  Let $f_B: [-1,1]\to [-1,1]$ denote the diffeomorphism given by $f_B(z) =  \Psi_B(1,z)$.

Note that since $X_B$ has positive slope whenever both $t$ and $z$ lie in $(-1+\epsilon, 1-\epsilon)$, $f_B(z)>z$ whenever $z$ lies in $(-1+\epsilon, 1-\epsilon)$. Indeed, by choosing the slope of $X_B$ to be great enough, we may guarantee that $f_B(-1/2)\ge 1/2$. Rechoose $X_B$ as necessary so that $X_B$ satisfies (1)--(3) and also:\\

(4) $f_B(-1/2)\ge 1/2$.\\

Now choose $X_D$ along the side $x=-1$ by setting $X_D(-1,y,z) = X_B(1,-y,z)$.

Notice that if $f :[-1,1]\to [-1,1]$ is any orientation preserving diffeomorphism, then there is a smooth  flow $\Psi=\Psi(x,z)$ on $[-1,1]\times [-1,1]$ such that $\Psi(1,z) = f(z)$; simply set
$$\Psi(x,z)=\frac{1-t(x)}{2} z + \frac{1+t(x)}{2} f(z)$$ 
 where $t(x)$ is a smooth function of $x$ that is $-1$ for $x\in[-1,-1+\epsilon]$, $1$ for $x\in [1-\epsilon, 1]$, and has positive derivative for all other $x$.
Thus, by specifying a diffeomorphism $f:[-1,1]\to [-1,1]$, we  specify a family of smooth flows $\Psi(x,z)$ and corresponding  smooth unit tangent vector fields $X$.

We take advantage of this to define  the vector field $X_A$ along the side $y = -1$. Let  $f_A: [-1,1]\to [-1,1]$ be the diffeomorphism given by
$$f_A^{-1}(z)=u\circ f_B \circ f_B\circ f_B(z),$$ where  
 $u: [-1,1]\to [-1,1]$ is a diffeomorphism which is the identity on $z\in [-1,-1+\epsilon]\cup [1-\epsilon,1]$ and strictly increasing  elsewhere.
 Let $\Psi_A$ be a smooth flow on $[-1,1]\times [-1,1]$
 which satisfies the following:
\begin{enumerate}
\item $\Psi_A$  has negative  slope whenever $z\in (-1+\epsilon, 1-\epsilon)$,
\item $\Psi_A$ has unit tangent vector field given by $\partial/\partial x$ when $z\in [-1,-1+\epsilon]\cup [1-\epsilon,1]$
\item $ \Psi_A(1,z) = f_A(z)$.
 \end{enumerate}
 Let $X_A$ be the smooth unit tangent vector field to $\Psi_A$.
 Recall that $\mathcal F$ is horizontal, and hence is dominated by $X_A$, along $A$.

Similarly, we use a diffeomorphism $f_C:[-1,1]\to [-1,1]$ to define a smooth vector field $X_C$ along the side $y=1$.
Let $f_C:[-1,1]\to [-1,1]$ be a diffeomorphism which satisfies:
\begin{enumerate}
\item $f_C(z)=1/2 f_A(2z)$ when $|z| \le 1/2$,
\item $f_C(z) = z$ when  $z\in [-1,-1+\epsilon]\cup [1-\epsilon,1]$, and
\item $f_C(z)<f_B(z)$ when $ 1/2 < |z| < 1-\epsilon$.
\end{enumerate}
Since  $f_A(z)=z$ whenever $|z|\ge 1-\epsilon$, $f_C(z)=z$ for $(1-\epsilon)/2 \le |z|\le 1/2$. Therefore, $f_C[-1/2,1/2]=[-1/2,1/2]$. Since $f_B(-1/2)>1/2$, it follows that $f_C(z)\le  f_B(z)$ for all $z$, with equality only when $|z|\ge 1-\epsilon$. So $f_B\circ f_C\circ f_B(z)<f_B\circ f_B\circ f_B (z)$, and hence 
$$f_A\circ f_B\circ f_C\circ f_B(z)< z,$$ 
for all $z \in (-1+\epsilon, 1-\epsilon).$

So  the diffeomorphism $f_A\circ f_B\circ f_C\circ f_B : [-1,1]\to [-1,1]$ is strictly decreasing on $(-1+\epsilon,1-\epsilon)$ and hence Proposition~\ref{cylinderextend}  applies.
\end{proof}

\begin{cor} \label{smoothisbracketed}  Let $\mathcal F$ be a taut, oriented, $C^2$-foliation of a closed oriented 3-manifold. Suppose $\mathcal F$ is not the product foliation $S^1\times S^2$. Then $\mathcal F$ is bracketed.
\end{cor}

\begin{proof} Let $\Phi$ be a volume preserving flow transverse to $\mathcal F$.
As noted in \cite{ET}, either $\mathcal F$ has a lot of nontrivial linear holonomy or it is $C^0$ close, and hence $\Phi$-close, to a foliation $\mathcal F'$ which has a lot of nontrivial linear holonomy.  In other words, there exists a disjoint union $V$ of  attracting neighborhoods such that one of $\mathcal F$ or $\mathcal F'$ is $V$-transitive.  In either case it follows that $\mathcal F$ is bracketed.
\end{proof}

\section{$L$-bracketed foliations}\label{L bracketed}

In this section, we introduce a new method for $\Phi$-approximating a foliation $\mathcal F$ by a pair of transitive confoliations, one positive and one negative. This method applies whenever there exists a link transverse to $\mathcal F$ which satisfies the condition  given in Definition~\ref{transtranv}. We also remark on some consequences yielding
3-manifolds containing weakly symplectically fillable contact structures.

Before stating Definition~\ref{transtranv}, we make some preliminary observations. 

\begin{prop} \label{flowextend}
Suppose $\mathcal F$ is a taut oriented codimension-1 foliation in $M$. Let $L$ be any link transverse to $\mathcal F$. Then there is a choice of metric on $M$ and volume preserving flow $\Phi$ everywhere transverse to $\mathcal F$ such that $L$ is contained as a union of closed orbits of $\Phi$. Moreover, given a choice of regular neighborhood $N(L) = \cup_i D_i\times S^1$ of $L$, the metric on $M$ and $\Phi$ can be chosen so that $\Phi$ is a trivial product $\{p\}\times S^1$, $p\in D_i$, on this regular neighborhood; in particular, $\Phi$ restricts to a flow on $\partial N(L)$.
\end{prop}

\begin{proof}
This follows immediately from the proof of Theorem~A1 found in \cite{Hass}.
\end{proof}

So any link transverse to $\mathcal F$ can be extended to a volume preserving flow $\Phi$ transverse to $\mathcal F$.  Alternatively, we may begin with a volume preserving flow $\Phi$ and let $L$ be a collection of closed orbits of $\Phi$. Without loss of generality, we will restrict attention to  the case that this flow $\Phi$ restricts to a flow on $\partial N(L)$ for some choice of regular neighborhood $N(L)$ of $L$. In either case, $\Phi$ determines a preferred, possibly non-compact, curve on each component of $\partial N(L)$:

\begin{lemma}\label{torusfoln}
Let $T$ be a framed torus and let $\Phi$ be a flow on $T$. Then either  $\Phi$ contains a simple closed curve  of some rational slope $m_{\Phi}$ or $\Phi$ is topologically conjugate either to a foliation by lines of  some irrational slope $m_{\Phi}$ or to a Denjoy blowup of a foliation by lines of  some irrational slope $m_{\Phi}$. In each case, the slope $m_{\Phi}$ is uniquely determined (by $\Phi$).
\end{lemma}

\begin{proof}  As long as no leaf of $\Phi$ has slope $1/0$, the framing determines a unique realization of  $\Phi$  as a suspension of some homeomorphism $f$ of  $S^1$  and the  Poincar\'e rotation number of $f$  determines the slope $m_{\Phi}$. See, for example, 4.3.1, 5.1.1 and 5.1.3 of  \cite{HH}.
\end{proof}

\begin{figure}[htbp] %  figure placement: here, top, bottom, or page
   \centering
   \includegraphics[width=3in]{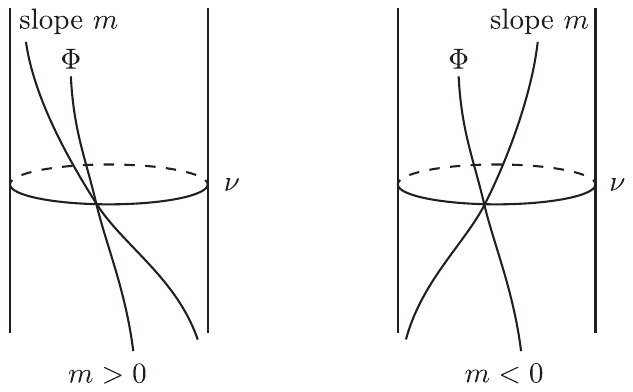}
   \caption{Slope convention on a component of $\partial N(L)$ in which slopes are designated as viewed from inside $N(L)$. }
   \label{slopeconvention}
\end{figure}

Denote this preferred isotopy class of curves, represented by either a simple closed leaf or an immersed $\mathbb R$, by $m_{\Phi}^T$.  We are interested in  the case that this torus $T$ is a component of $\partial N(L)$. In this case, there is also the isotopy class of the meridian, $\nu^T$ say, and $\nu^T\ne m_{\Phi}^T$. We shall call an isotopy class of a nontrivial curve $C$ in $T$ {\sl positive} if it has positive slope with respect to $\langle \nu, m_{\Phi}\rangle$ when viewed from inside $N$. Similarly, we shall call an isotopy class of curves {\sl negative} if it has negative  slope with respect to $\langle \nu, m_{\Phi}\rangle$ when viewed from inside $N$. This convention is illustrated in  Figure~\ref{slopeconvention}.

Note that if $\mathcal F$ is an oriented codimension-1 foliation which intersects a torus $T$ transversely, then $\mathcal F\cap T$ is a flow on $T$, and  $m_{\mathcal F}^T$ denotes  the preferred isotopy class of this flow.

\begin{definition} 
Define a triple $(M,\mathcal F, \Phi)$ to be {\sl coherent} if the foliation $\mathcal F$ is taut and oriented, the flow $\Phi$ is volume preserving and positively transverse to $\mathcal F$, and the boundary of M, if nonempty, is a union of flow lines.
\end{definition}

Let $L$ be a link in $M$.  A foliation $\mathcal F$ is {\sl $L$-taut} if $\mathcal F\pitchfork L$, and $\mathcal F$ is $L$-transitive, that is, each leaf of $\mathcal F$ has nonempty intersection with $L$. Similarly, if $\partial M\ne \emptyset$, then $\mathcal F$ is {\sl $\partial M$-taut} if $\mathcal F$ intersects $\partial M$ transversely, with no Reeb annuli, and $\mathcal F$ is $\partial M$-transitive, that is, each leaf of $\mathcal F$ has nonempty intersection with $\partial M$.

Recall that a foliation $\mathcal F_0$ is said to {\sl realize slope} $m$ along a framed torus boundary component $T$ if $\partial \mathcal F_0\cap T$ consists of parallel curves, not necessarily compact, of slope $x$. When $\mathcal F_0$ is oriented, these curves $\partial \mathcal F_0\cap T$ are necessarily consistently oriented. Notice that the condition that a foliation $\mathcal F_0$ be $\partial M_0$-taut is weaker than the condition that $\partial F_0$ realizes slope $m_{\mathcal F_0}^T$ for each component $T$ of $\partial M_0$; in other words, nontrivial holonomy is possible for $\mathcal F_0\cap \partial M_0$.

\begin{definition}\label{transtranv}
Suppose $\mathcal F$ is a taut oriented codimension-1 foliation in $M$.  Let $L$ be a link in $M$ which is transverse to $\mathcal F$ and let $M_0$ equal $M\setminus \text{int}N(L)$. Let $\mathcal F_0$ denote the restriction of $\mathcal F$ to $M_0$. The foliation $\mathcal F$ is  {\sl $L$-bracketed} if, for some choice of metric on $M_0$, there is a volume preserving flow $\Phi_0$ on $M_0$ such that $(M_0,\mathcal F_0,\Phi_0)$ is coherent and the following property is satisfied:

$M_0$ contains a pair of foliations $\mathcal F_{\pm}$ such that 
\begin{enumerate}

\item $(M_0,\mathcal F_+,\Phi_0)$ and  $(M_0,\mathcal F_-,\Phi_0)$ are coherent,

\item $\mathcal F_{\pm}$ are $\partial M_0$-taut, and

\item for each component $T$ of $\partial M_0$, $m_{\mathcal F_-}^T$ is negative and $m_{\mathcal F_+}^T$ is positive with respect to $\langle \nu^T, \Phi|_T\rangle$, where $\nu^T$ is the meridian slope of component $T$ (and hence is the slope of $\mathcal F_0\cap T$).
\end{enumerate}
To make the flow explicit, we also say that $M$ contains an {\sl $(\mathcal F,\Phi)$-transitive} link $L$, where $\Phi$ is the flow $\Phi_0$ blown down to a (volume preserving) flow on $M$.
\end{definition}

The notion of $L$-bracketed is a special case of a bracketed foliation $\mathcal F$.  The decomposition is given by setting $V=V'=N(L)$, $W=W'=M_0$.  There is a canonical choice of positive or negative contact structure on $N(L)$ which is given by perturbing the meridional disks, and then the requirement is to find foliations $\mathcal F_W=\mathcal F_+$ and $\mathcal F_{W'}=\mathcal F_-$ on $M_0$.

\begin{thm} \label{transitivemain}
Suppose $\mathcal F$ is a taut oriented codimension-1 foliation in $M$ and that $\mathcal F$ is $L$-bracketed for some link $L$. Then $\mathcal F$ can be $\Phi$-approximated by a pair of smooth contact structures $\xi_{\pm}$, one positive and one negative.  
  These contact structures $(M,\xi_+)$ and $(-M,\xi_-)$ are weakly symplectically fillable and universally tight. 
\end{thm}

\begin{proof}
Set $V= N(L)$, and let $W$ denote the closure of the complement of $V$.  Define a contact structure $\xi_0$ on $N(L)$ so that each component of $L$ is a transverse knot and each component of $N(L)$ is a standard positive contact neighborhood of its core.  Choosing the rate of rotation of the contact planes along each meridional disk to be small guarantees that the characteristic foliation of $\xi_0$ along $\partial N(L)$ is close to the meridian. It follows that we may choose  $\xi_0$   so that it is  strictly dominated by $\mathcal F_+$. Apply Theorem~\ref{main1} to obtain $\xi_+$. Similarly, each component of $N(L)$ can be modeled using the standard negative radial model, and Theorem~\ref{main1} can be applied to obtain $\xi_-$. 

Since $\xi_{\pm}$ are both positively transverse to the volume preserving flow $\Phi$, they are weakly symplectically fillable.
\end{proof}

Since transitive links are somewhat mysterious, it is natural to ask:

\begin{question} Given a foliation $\mathcal F$, does there exist  a link $L$ for which $\mathcal F$ is $L$-bracketed?
\end{question}

It is not clear how the answer to this might change if the link is required to be connected.

\begin{question} Given a foliation $\mathcal F$, does there exist  a knot $K$ for which $\mathcal F$ is $K$-bracketed?
\end{question}

\begin{example}\label{s1timess2}  The product foliation $\mathcal F$ on $S^1 \times S^2$ is an example of a foliation which is not $L$-bracketed for any link $L$.  The existence of such a link would imply, by Theorem~\ref{transitivemain}, that $\mathcal F$ can be approximated by a tight contact structure $\xi$.  This would imply that the underlying 2-plane bundles of $\mathcal F$ and $\xi$ are equivalent.  The Euler class of the foliation, $e(\mathcal F)$, evaluated on a spherical leaf $S^2$ of $\mathcal F$ equals 2.  On the other hand, this $S^2$ is homotopic to a convex surface $S'$, \cite{Gi1}, in a tight contact structure.  It follows that $S'$ has a connected dividing set, \cite{Gi2}, and that $e(\xi)$ vanishes on it. 
\end{example}

A taut foliation $\mathcal F$ is certainly $L$-taut for some link $L$. In fact, it is $L$-taut for some knot $L$. Moreover, as noted above, for some choice of metric there is a volume preserving flow $\Phi$ transverse to $\mathcal F$ and containing $L$ as an orbit or union of closed orbits. And often, although not necessarily, foliations on $M_0$ that are $\Phi_0$-close to $\mathcal F$ will also be $L$-taut. So a key question is the existence of a pair of $\Phi_0$-close foliations $\mathcal F_{\pm}$ in $M_0$ such that $m_{\mathcal F_-}^T< \text{slope}\, \partial \mathcal F_0|_T< m_{\mathcal F_+}^T$ for each component $T$.

\begin{example} Consider the case that $M_0$ is any compact orientable manifold with boundary a nonempty union of  $b$ tori. Suppose that $\Phi_0$ is a volume preserving flow which is tangent to $\partial M_0$ and  that $B$ is a  transversely oriented branched surface transverse to $\Phi_0$. 

If $B$ fully carries a set of foliations which are $\partial M_0$-taut and realize a nonempty open set $J$ of boundary slopes (if $b=1$) or multi-slopes (if $b\ge 2$), then Dehn-filling $M_0$ along any rational slope or multi-slope in $J$ results in a foliation which is $L$-bracketed, where $L$ is the link which is the core of the Dehn filling. Examples of such foliations can be found in the papers \cite{DL,DR,g1,g2,g3,G3,KRo,KR2,Li2,LR,R,R1,R2}.

One can ask whether the foliations constructed by Dehn filling more than one torus can be $K$-bracketed for some knot $K$. 

In \cite{G1, G2,G3} Gabai constructs foliations in closed manifolds $M$ with $H_2(M)\ne 0$. These foliations are fully carried by finite depth branched surfaces. One can ask whether such foliations are $L$-bracketed for some link $L$.
\end{example}
 
Finally, we note that the proof of Theorem~\ref{transitivemain} doesn't  actually require the existence of the foliation $\mathcal F_0$. More precisely, we have the following.

\begin{definition}\label{tranvlink}
Suppose $(M,\xi)$ is a contact 3-manifold.  Let $L$ be a transverse link in $(M,\xi)$. Let $M_0$ equal $M\setminus \text{int}N(L)$. The contact structure $\xi$ is  {\sl $L$-bracketed} if, for some choice of metric on $M_0$, there is a volume preserving flow $\Phi_0$ on $M_0$, tangent to $\partial M_0$, such that the following property is satisfied:

$M_0$ contains a pair of foliations $\mathcal F_{\pm}$ such that 
\begin{enumerate}

\item $(M_0,\mathcal F_+,\Phi_0)$ and  $(M_0,\mathcal F_-,\Phi_0)$ are coherent,

\item $\mathcal F_{\pm}$ are $\partial M_0$-taut, and

\item for each component $T$ of $\partial M_0$, $m_{\mathcal F_-}^T$ is negative and $m_{\mathcal F_+}^T$ is positive with respect to $\langle \nu^T, \Phi|_T\rangle$, where $\nu^T$ is the meridian slope of component $T$.\end{enumerate}

\end{definition}

\begin{thm} \label{transitivemaincontact}
Suppose  $(M,\xi)$ is a contact 3-manifold and that $\xi$ is $L$-bracketed for some transverse link $L$. Then $\xi$ can be $\Phi$-approximated by a pair of smooth contact structures $\xi_{\pm}$, one positive and one negative. These contact structures are necessarily weakly symplectically fillable and universally tight.
\end{thm}

\section{Open book decompositions}\label{Open book}

An interesting class of $L$-bracketed foliations is obtained by considering the special case that $L$ is a fibered link in $M$ and $\mathcal F$ is transverse to a flow $\Phi$  obtained by surgery from a volume preserving suspension flow of the corresponding fibre bundle complement of $L$. In this case $L$ forms the binding of an open book decomposition $(S,h)$ of $M$ and the contact structure $\xi_{(S,h)}$ compatible with $(S,h)$ is $\Phi$-close to $\mathcal F$. 

For completeness, we begin with some standard definitions. Since we are relating ideas from the world of codimension-1 foliations and the world of contact structures, we will also provide some translations between the terminologies of these two worlds. The main results of this section appear in Subsection~\ref{OBresults}.

\subsection{Open book decompositions}

Let $S$ be a compact surface with nonempty boundary. A pair $(S,h)$, where $h$ is a homeomorphism that restricts to the identity map on $\partial S$, determines a closed 3-manifold $M=S\times [0,1]/\approx$ where the equivalence relation $\approx$ identifies $(x,1)\approx(h(x),0)$ for all $x\in S$ and $(x,s)\approx(x,t)$ for all $x\in \partial S$ and $s,t\in [0,1]$. The singular fibration with pages $S\times\{t\}$ is called the {\sl open book determined by the data $(S,h)$,} and we write $M=(S,h)$.

\subsection{Surface bundles over $S^1$}

Corresponding to  an open book decomposition of $M$ is a description of $M$ as a Dehn surgery along the binding $L=\cup l_i$ by meridional multislope $(\nu_1,...\nu_b)$.
Conversely, corresponding to such a Dehn filling description of $M$, we have a corresponding open book description of $M$. Since many existing constructions of foliations are
described from the Dehn surgery perspective, it is useful to consider this correspondence more carefully.

Let $M_0$ denote the compact complement of $L$; so
$M_0= M\setminus \text{int} N(L)$, where $N(L)$ is a regular neighbourhood of $L$, and $M_0$  is homeomorphic to $S\times [0,1]/h$.

Notice that if $S$ is a disk, then necessarily  $h$ is isotopic rel boundary to the identity map. If $S$ is an annulus, then $h$ is isotopic rel boundary to  some power of the Dehn twist about the core of $S$. Otherwise, $S$ is hyperbolic. We therefore lose little by restricting attention to the case that $S$ is hyperbolic and will now do so.

Recall Thurston's classification of surface automorphisms. 

\begin{thm}\label{Thurston} \cite{Th, CB, FLP} Let $S$ be an oriented hyperbolic surface with geodesic boundary, and let $h\in \text{Homeo}(S,\partial S)$. Then $h$ is freely isotopic to either 

(1) a pseudo-Anosov homeomorphism $\theta$,

(2) a periodic homeomorphism $\theta$, in which case there is a hyperbolic metric for which $S$ has geodesic boundary and such that $\theta$ is an isometry of $S$, or

(3) a reducible homeomorphism $h'$ that fixes, setwise, a maximal collection of  disjoint simple closed geodesic curves $\{C_j\}$ in $S$.
\end{thm}

Recall that a pseudo-Anosov homeomorphism has finitely many prong singularities and is smooth and hyperbolic elsewhere \cite{FLP}. To avoid overlap in the cases, we refer to a map as reducible only if it is not periodic. Since we will be considering homeomorphisms $h$ in the context of open books $(S,h)$, we will be considering only homeomorphisms $h$ which fix $\partial S$ pointwise. Therefore, given a reducible map, splitting $S$ along $\cup_jC_j$ gives a collection of surfaces $S_1,\dots, S_n \subset S$ with geodesic boundary that are fixed by $h'$.   Maximality of $\{C_j\}$ implies that applying Thurston's classification theorem to each $h'|_{S_i}\in \text{Homeo}(S_i, \partial S_i)$ produces either a pseudo-Anosov or periodic representative. So we may assume that $h'$ is either periodic or pseudo-Anosov away from some small neighborhood of the $C_i$.

\begin{definition} \label{Thurstonrep} Let $S$ be hyperbolic and $h\in \text{Homeo}(S,\partial S)$. If conclusion (1) or  conclusion (2) of Theorem~\ref{Thurston} is satisfied, call $\theta$  \textit{the Thurston representative} of $h$. If instead conclusion (3) holds, let $\theta: (S,\partial S)\to (S,\partial S)$ denote the piecewise continuous function uniquely determined by the following constraints:
\begin{itemize}

\item [(3.1)] $\theta$ restricted to each component of the complement of the union $\cup C_i$ is freely isotopic to the restriction of $h'$ to this component and is either periodic or pseudo-Anosov. 

\item [(3.2)] $\theta$ restricted to each simple closed geodesic $C_i$ is  freely isotopic to the restriction of $h'$ to $C_i$ and is a periodic isometry, and 

\end{itemize} Again, refer to $\theta$ as the {\sl Thurston representative} of $h$.
\end{definition}

Now consider again the open book decomposition $M=(S,h)$ and let $\theta$ denote the Thurston representative of  $h$. When $\theta$ is periodic or pseudo-Anosov, the link complement $M_0 = M\setminus \text{int} N(L)$ is also homeomorphic  to the mapping torus $S\times [0,1]/(x,1)\sim (\theta(x),0)$ of $\theta$, and in the discussions which follow, we will typically view $M_0$ as the mapping torus of $\theta$. When $\theta$ is reducible and so only piecewise continuous, we will typically view $M_0$ as the union along essential tori of the mapping tori of the extension to $\cup S_i$ of the restriction of $\theta$ to  the complement of the union $\cup C_i$.
 
Let  $\Theta_0$  be the flow  obtained by integrating the vector field $\partial/\partial t$, where points of $M_0$ are given by $[(x,t)], x\in S, t\in [0,1]$.  We will refer to this flow as either the  {\it suspension flow of $\theta$} or the {\it Thurston flow  (associated to $h$)}.

Since $\theta$ is area preserving with respect to some metric on the fibre, $\Theta_0$ is volume preserving with respect to some choice of metric on $M_0$.

Notice that the suspension flow $\Theta_0$ is pseudo-Anosov (respectively, periodic) when $\theta$ is pseudo-Anosov (respectively, periodic). In particular, when $\theta$ is periodic, all orbits of $\Theta_0$ are closed.  When $\theta$ is pseudo-Anosov, there are an even number of alternately attracting and repelling closed orbits of $\Theta_0|{\partial N(l_i)}$ for each $i$. When $\theta$ is periodic or pseudo-Anosov, this flow is continuous and the orbits are smoothly embedded.  When $\theta$ is reducible, this flow is not continuous, but orbits of the flow are smoothly embedded. 

Since many closed manifolds, together with corresponding open book decompositions, can be realized by Dehn filling $M_0$, it is useful to have canonical framings on the boundary components of $M_0$ which are defined independently from $M$. As in \cite{R2}, we will use the Thurston flow $\Theta_0$ to define these canonical coordinate systems on $\partial M_0$.

\subsection{The Thurston flow framing on surface bundles over $S^1$}

\begin{definition} \label{coordinatesgalore} Let $\partial_iM_0$ denote the $i$-th boundary component of $M_0$. Choose an oriented identification $\partial_i M_0 \sim \mathbb R^2/\mathbb Z^2$ by choosing oriented curves $\lambda_i$ and $\mu_i$, so that $\lambda_i$ has slope $0$ and $\mu_i$ has slope $\infty$, as follows. Let $\lambda_i = \partial (S\times \{0\})$, with orientation induced by the orientation on $S$. Let $\gamma_i$ be a closed orbit of the flow $\Theta_0$ restricted to $\partial_iM_0$. Choose $\mu_i$ to be an oriented  simple closed curve which has algebraic intersection number $< \lambda_i,\mu_i> = 1$  and which minimizes the geometric intersection number $|\gamma_i\cap\mu_i|$. This choice is unique except in the case that the geometric intersection number $|\gamma_i\cap\lambda_i| = 2$. In this case we choose $\mu_i$ so that $\gamma_i$ has slope $+2$.  Call the resulting framing the \textit{Thurston flow framing} on $\partial_iM_0$.  Slopes expressed in terms of the  flow framing will be said to be given in {\sl Thurston  flow coordinates}. This was originally called the {\sl natural framing} or {\sl natural coordinates} by Roberts in \cite{R2}.
In this section,  we are beginning with a fixed fibering and hence the associated Thurston flow coordinates are well defined. In general, different choices of fibering can lead to nonisotopic  closed orbits $\gamma_i$, and hence to different Thurston flow coordinates. Notice also that in these coordinates, the slope of $\gamma_i$ always satisfies  $$1/ (\mbox{slope  } \gamma_i)  \in (-1/2,1/2].$$
 
\end{definition}

Now let's consider the relationship between the flow framing and the fractional Dehn twist coefficient defined by Honda, Kazez, and Mati\'c in \cite{HKMRV1}. First recall the definition of fractional Dehn twist coefficient.

\begin{definition} \cite{HKMRV1}
Fixing a component $C_i$ of $\partial S$ and  restricting the flow $\Theta_0$ to the component of  $\partial N(L)$ corresponding to $C_i$, $\Theta_0$  necessarily has periodic orbits.  Let $\gamma_i$ be one such, and write $$\gamma_i = p_i\lambda_i+q_i\nu_i$$ where  $\lambda_i = C_i$, $\nu_i$ is the meridian (oriented so that $<\lambda_i,\nu_i>=1$), and  $p_i$ and $q_i$ are relatively prime integers with $q_i>0$.  The {\sl fractional Dehn twist coefficient} of $h$ with respect to the component $C_i$ of $\partial S$ is given by $$c_i(h)=p_i/q_i\, .$$  In particular, when $\gamma_i=\nu_i, (p_i,q_i)=(0,1)$ and $c_i(h)=0.$
\end{definition}

Recall that  $M$ is obtained from $M_0$ by $(\nu_1,...,\nu_b)$ Dehn filling along the boundary components of $M_0$.  Beginning with the open book decomposition  of $M$ and the associated  fractional Dehn twist coefficients $c_i=p_i/q_i, 1\le i\le b,$ (as above, expressed in $(\lambda_i,\nu_i)$ coordinates), it is sometimes  useful to express the slopes of $\nu_i$ and $\gamma_i$ in terms of the Thurston flow coordinates, $(\lambda_i,\mu_i)$. We now describe how to do this.  In flow coordinates, $\lambda_i$ has slope $0$. Since $|\lambda_i\cap\nu_i| = 1$, it follows that, in flow coordinates, $\nu_i$ has slope $1/k_i$ for some integer $k_i$.  As noted in \cite{HKM2}, the integer $k_i$ is uniquely determined by the fractional Dehn twist coefficient $c_i(h)$ for each $i, 1\le i \le b$.   This relationship can be very simply stated:

\begin{prop} \label{translation1}{\bf(Coordinate translation I)} Let $c_i = c_i(h)$ and let $n_i$ be the integer determined by the condition
$$c_i\in (n_i-1/2,n_i+1/2].$$ 
In other words, $$n_i=\lceil c_i -1/2 \rceil,$$  the integer nearest to $c_i$, with ties in the case $c_i\in \mathbb Z+1/2$ broken by rounding down.

Then $k_i = -n_i$ and so $\nu_i$ has slope $-1/n_i$. Moreover, $\gamma_i$ has slope  $1/(c_i-n_i)$ and so
\begin{enumerate}
\item if $c_i = n_i$, then  $\gamma_i=\mu_i$ has slope $1/0$;

\item if $c_i > n_i$, then $\gamma_i$ has positive slope; and

\item if $c_i<n_i$, then $\gamma_i$ has negative slope.
\end{enumerate}
\end{prop}

\begin{proof}
The meridian $\nu_i$ has slope $1/k_i$ and so $\nu_i = \mu_i + k_i\lambda_i$. 
So $$\gamma_i = p_i\lambda_i+ q_i\nu_i = p_i\lambda_i + q_i(\mu_i + k_i\lambda_i) = q_i\mu_i +(p_i+k_iq_i)\lambda_i $$
has slope   ${q_i}/(p_i+k_iq_i) = {1}/(c_i+k_i)$ and $|\gamma_i\cap\mu_i| =  |p_i+k_iq_i|$. By definition of flow coordinates, $k_i$ is chosen to minimize $|\gamma_i\cap\mu_i| =  |p_i+k_iq_i|$, and hence to minimize $|c_i + k_i|$. There is a unique such minimizing $k_i$ unless $c_i\in \mathbb Z + 1/2$. In this case, $k_i$ is chosen so that $\gamma_i$ has slope $2/1$; namely, so that  $c_i = -k_i + 1/2$. So $k_i$ is the unique integer satisfying $c_i\in (-k_i-1/2,-k_i + 1/2]$.
\end{proof}

Conversely, given the fibered 3-manifold $M_0$ and meridional  Dehn filling slopes $\nu_i, 1\le i\le b,$ in terms of the Thurston flow coordinates $(\lambda_i,\mu_i)$, it is often useful to express the slopes $\nu_i$ and $\gamma_i$ in terms of the associated open book coordinates $(\lambda_i,\nu_i)$.  We have the following. 

\begin{prop} \label{translation2}{\bf(Coordinate translation (II)}
Suppose $M_0$ is fibered, with the boundary components $\partial_i M_0$ given the Thurston flow framing $(\lambda_i,\mu_i)$ for each $i$. As above, let $\nu_i$ be a meridional slope and let $\gamma_i$ be a closed orbit of the Thurston flow on $\partial_i M_0$.  In terms of the Thurston framing, $\nu_i=-1/n_i$ and $\gamma_i=r_i/s_i$, for some integers $n_i,r_i$ and $s_i$. Again as above, let $M$ be the manifold obtained by $(\nu_1,...,\nu_b)$ filling $M_0$ and let $(S,h)$ be the open book decomposition of $M$ determined by the fibering of $M_0$. Then, in terms of the open book framing $(\lambda_i,\nu_i)$ on $\partial_i M_0 = \partial N(l_i)$, 
$$\nu_i=1/0, \,\,\, \mu_i = 1/n_i \,\,\,  \mbox{  and  } \gamma_i = r_i/(n_i r_i + s_i).$$ In particular, the fractional Dehn twist  coefficient along  or $S \cap \partial_iM_0$ is given by 
$$c_i(h) = n_i+s_i/r_i, $$where, as noted in
Definition~\ref{coordinatesgalore}, $s_i/r_i\in (-1/2,1/2]$.
\end{prop}

\begin{proof}
To eliminate the ugliness of subscripts, focus on a particular boundary component, $\partial_i M_0$, and drop all reference to $i$.

Two right handed framings are related by a unique transformation  in $SL_2(\mathbb Z)$.  Notice that 
\[
   A=
  \left[ {\begin{array}{cc}
   1 & n \\
   0 & 1 \\
  \end{array} } \right]
\]
is the  element in $SL_2(\mathbb Z)$ which maps the pair  $({1 \choose 0},{-n \choose 1})$ to the pair $({1 \choose 0},{0 \choose 1})$.

Hence, with the correspondence ${a\choose b} \mapsto b/a$,  $A$ describes the translation from $(\lambda,\mu)$ coordinates to $(\lambda,\nu)$ coordinates.  The slope computations follow immediately.
\end{proof}

\subsection{Dehn filling the Thurston flow}  Dehn surgery on Anosov flows is defined by Goodman in  \cite{Go}. (See also \cite{Fried}.) This definition generalizes naturally to the setting of pseudo-Anosov flows and permits us to consider the effect of Dehn filling of  Thurston flows. We define Dehn filling of a  Thurston flow as follows.  Let $Y$ be any closed 3-manifold obtained by Dehn filling $M_0$.  For each $i, 1\le i\le b$, let $X_i$ denote the solid torus in $Y$ bounded by $\partial_i M_0$ and let $\kappa_i$ denote the core of $X_i$.  As long as the surgery coefficient along $\partial_i M_0$ is not $\gamma_i$, it is possible to blow down $X_i$ to its core and obtain a flow $\Theta$ defined on $Y$. Notice that the cores $\kappa_i$ are closed orbits of $\Theta$. Also, either  $\Theta_0$ is periodic in a neighborhood of $\partial_i M_0$ and therefore  so is $\Theta$ in a neighborhood of $\kappa_i$, or else $\Theta_0$ is pseudo-Anosov (with possibly a single prong pair along $\kappa_i$) in a neighborhood of $\partial_i M_0$ and therefore  so is $\Theta$ in a neighborhood of $\kappa_i$. We shall refer to this flow $\Theta$ as the {\sl surgered Thurston flow}. Notice that since $\Theta_0$ is volume preserving for some metric, so is $\Theta$.
 
\begin{lemma} Suppose $M$ has an open book decomposition $(S,h)$ with binding $L$ and corresponding Dehn filling description $M=M_0(\nu_1,\cdots,\nu_b)$.  Let $\Theta_0$ denote the Thurston flow on the complement of $L$. Then $\Theta_0$ extends to a flow $\Theta$ on $M$ if and only if all fractional Dehn twist coefficients $c_i$ are nonzero.
\end{lemma}

\begin{proof}
Notice that $$c_i(h)=0 \iff p_i = 0, q_i = 1 \iff \gamma_i = \nu_i \implies \nu_i = \mu_i = \gamma_i.$$
In particular, if at  every component $C_i$ of $\partial S$ the fractional Dehn twist coefficient $c_i(h)\ne 0$, then $\nu_i\ne\gamma_i$ and so it is possible to blow down the flow $\Theta_0$ to the surgered Thurston flow $\Theta$ on $M$. Otherwise it is not. 
\end{proof}

Notice that the binding $L=\cup l_i$ inherits orientations both from the open book structure and from the flow $\Theta$. 

\begin{figure}[htbp] %  figure placement: here, top, bottom, or page
   \centering
   \includegraphics[width=4.2in]{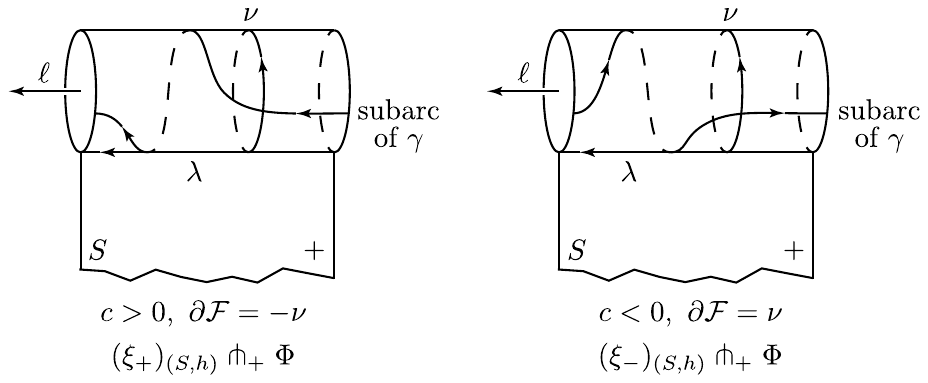} 
   \caption{ Curves on a boundary component of $M_0$.}
   \label{postrans}
\end{figure}

\begin{lemma} \label{bindingorientation}
When $c_i>0$, these orientations agree on $l_i$. When $c_i<0$, these orientations do not agree on $l_i$.
\end{lemma}

\begin{proof}  The orientation of $\Theta$ restricted to the binding is determined by the sign of the slope of $\gamma_i$ as expressed in $(\lambda_i,\nu_i)$ coordinates. This is illustrated in Figure~\ref{postrans}.
\end{proof}

\subsection{Contact structures supported by an open book}

In \cite{Gi}, Giroux defined the notion of contact structure supported by an open book: a positive (respectively, negative) contact structure $\xi$ on $M$  is {\sl supported by, or compatible with, an open book decomposition $(S,h)$ of $M$} if $\xi$ can be isotoped through contact structures so that there is a contact 1-form $\alpha$ for $\xi$ such that 
\begin{enumerate}
\item $d\alpha$ is a positive area form on each page $S_t$ of the open book and
\item $\alpha>0$ (respectively, $\alpha<0$) on the binding $l$.
\end{enumerate}

\begin{prop} (\cite{Gi})
Two contact structures supported by the same open book are contact isotopic.
\end{prop}

We may therefore abuse language and refer to {\it the} contact structure compatible with the open book decomposition $(S,h)$.

\begin{notation} Let $(\xi_+)_{(S,h)}$ denote the positive contact structure compatible with the open book decomposition $(S,h)$. Let $(\xi_-)_{(S,h)}$ denote the negative contact structure compatible with the open book decomposition $(S,h)$.
\end{notation}

\begin{lemma}\label{posnegc} Suppose all fractional Dehn twist coefficients are nonzero and let $\Theta$ denote the surgered Thurston flow.
\begin{enumerate}
\item $(\xi_+)_{(S,h)}$ is positively transverse to $\Theta$ if and only if all fractional Dehn twist coefficients are positive.
\item $(\xi_-)_{(S,h)}$ is positively transverse to $\Theta$ if and only if all fractional Dehn twist coefficients are negative.
\end{enumerate}
\end{lemma}

\begin{proof} This follows immediately from Lemma~\ref{bindingorientation}.
\end{proof}

Any contact structure $(M,\xi)$ is supported by infinitely many open book decompositions. Honda, Kazez, and Mati\'c proved in \cite{HKM2} that if there is a compatible open book decomposition with nonpositive fractional Dehn twist coefficient $c_i$ for some $i$, then necessarily $\xi$ is overtwisted. In fact, they show the following.

\begin{thm} [Theorem 1.1 of \cite{HKMRV1}]
A contact structure $(M,\xi)$ is tight if and only if all of its compatible open book decompositions $(S,h)$ have fractional Dehn twist coefficients  $c_i \ge 0$ for $1\le i\le |\partial S|$.
\end{thm}

Let $S'$ be the smallest invariant subsurface of $S$ for the Thurston representative of $h$; thus $S'=S$ if and only if $h$ is not reducible. If $c_i=0$ for some common boundary component of $S$ and $S'$, there are two possibilities for $h|S'$.  The first is that it is periodic, and it follows that it is equal to the identity map.  The second possibility is that $h|S'$ is pseudo-Anosov, but this immediately implies $\xi$ is overtwisted \cite{HKMRV1}.  Since our primary focus is studying tightness for pseudo-Anosov maps it is enough to consider open book decompositions for which all fractional Dehn twist coefficients $c_i$ are positive (or, in the case that the contact structure is negative, to a consideration of open book decompositions for which all fractional Dehn twist coefficients $c_i$ are negative).

\subsection{Foliations compatible with an open book decomposition}\label{OBresults}

\begin{definition}
Let $\mathcal F$ be an oriented foliation of $M$. Let $(S,h)$ be an open book decomposition of $M$, with binding $L$. Let  $\Phi$ be the surgered Thurston flow associated to $h$. If $\mathcal F$ is everywhere transverse to $\Phi$, then we say that $\mathcal F$ is {\sl compatible with the open book} $(S,h)$.
\end{definition}

\begin{prop} \label{oppfoln} Let $(S,h)$ be an open book decomposition of $M$, with binding $L$, and surgered Thurston flow $\Phi$. Set $M_0=M\setminus \text{int} N(L)$ and $\Phi_0$ to be the Thurston flow associated to $h$.
\begin{enumerate}

\item Suppose all fractional Dehn twist coefficients are positive. If there exists a foliation $\mathcal F_-$ in $M_0$ which is $L$-taut, is transverse to $\Phi_0$ and satisfies $m_{\mathcal F_-}^T < 0$ with respect to $\langle \nu^T,\Phi|_T \rangle$ for all components $T$ of $\partial M_0$, then $  (M,(\xi_+)_{(S,h)}) $ is weakly symplectically fillable   and universally tight.

\item Suppose all fractional Dehn twist coefficients are negative. If there exists a foliation $\mathcal F_+$ in $M_0$  which is $L$-taut, is  transverse to $\Phi_0$ and satisfies $m_{\mathcal F_+}^T > 0$ with respect to $\langle \nu^T,\Phi|_T \rangle$ for all components $T$ of $\partial M_0$, then $  (-M,(\xi_-)_{(S,h)}) $ is weakly symplectically fillable   and universally tight.

\end{enumerate}
\end{prop}

\begin{proof} Consider case (1). Since $(\xi_+)_{(S,h)}$ is transverse to $\Phi$, it suffices to show that there is a negative contact structure $\xi_-$ which is  transverse to $\Phi$. This follows immediately from the assumptions   by Theorem~\ref{main1}.

Case (2) follows similarly.
\end{proof}

Set $\mathcal G_0$ to be the fibration $S\times [0,1]/h$ in the complement of $L$. Note that  $\mathcal G_0\cap T$ has either positive or negative slope with respect to $\langle \nu^T,\Phi|_T \rangle$, on each component $T$ of $\partial N(L)$. For completeness, we note the following:
 
\begin{lemma}
Let $T_i$ be a component of  of $\partial N(L)$. Then $\mathcal G_0\cap T_i$ has  positive (respectively, negative) slope with respect to $\langle \nu^{T_i},\Phi|_{T_i}\rangle$ if the corresponding fractional Dehn twist $c_i$ is positive (respectively, negative).
\end{lemma}
 
\begin{proof} Consider the relative slope values of $\gamma, \lambda$, and $\nu$ on any boundary component $T_i$ of $\partial M_0$. This is captured in Figure~\ref{postrans}. Notice that $\lambda$ represents the slope of $\partial \mathcal G_0$ and $\gamma$ represents the slope of $\Phi$. Recalling  the slope convention,  illustrated in Figure~\ref{slopeconvention}, we see that  $\mathcal G_0\cap T_i$ has  positive (respectively, negative) slope with respect to $\langle \nu^{T_i},\Phi|_{T_i}\rangle$ if the corresponding fractional Dehn twist $c_i$ is positive (respectively, negative).
\end{proof}
 
So $\mathcal G_0$ extends to a positive confoliation transverse to $\Phi$ when $c_i>0$ for all $i$ and to a negative confoliation transverse to $\Phi$ when $c_i<0$ for all $i$. In other words, and unsurprisingly, $\mathcal G_0$ as a foliation playing the role of $\mathcal F_{+}$  (respectively, $\mathcal F_-$) gives a second way of establishing the existence of $\xi_+$ (respectively $\xi_-$).

Finally, we use Proposition~\ref{translation1} to rephrase  Proposition~\ref{oppfoln} in terms of  Thurston flow coordinates. For simplicity of exposition, we restrict attention to the case that all fractional Dehn twist coefficients are positive. There is a symmetric statement in the case that all fractional Dehn twist coefficients are negative.

\begin{prop}\label{oppfolnTh}
Let $(S,h)$ be an open book decomposition of $M$, with binding $L$.  Suppose all fractional Dehn twist coefficients $c_i, 1\le i\le b,$ are positive. For each $i$, let $n_i$ be  the integer nearest to $c_i$, with ties in the case $c_i\in \mathbb Z+1/2$ broken by rounding down.
Set $M_0=M\setminus \text{int} N(L)$. Let $\Phi_0$ denote the Thurston flow associated to $h$ and let $\Phi$ denote the surgered Thurston flow on $M$.  

Suppose there exists a foliation $\mathcal F_-$ in $M_0$ which is $L$-taut, is transverse to $\Phi$, and such that, for each component $T_i$ of $\partial M_0$, $c_i$, $n_i$ and $$m_i = m^{T_i}_{\mathcal F_-}$$ satisfy one of the following:
\begin{enumerate}
\item $c_i=n_i$ and $m_i\in (-\infty,-\frac{1}{n_i})$, 
\item $c_i>n_i$ and $m_i\in (\frac{1}{c_i-n_i},\infty] \cup [-\infty,-\frac{1}{n_i})$, or
\item $c_i<n_i$ and $m_i\in (\frac{1}{c_i-n_i},-\frac{1}{n_i})$.
\end{enumerate}

Thus $m_{\mathcal F_-}^T<0$ with respect to $\langle \nu^T,\Phi|_T \rangle$ for all components $T$ of $\partial M_0$, and consequently, $(\xi_+)_{(S,h)}$ is weakly symplectically fillable   and universally tight.
\end{prop}

\begin{proof} The boundary slope of $\mathcal F_-$ on the $i^{th}$ boundary component lies between $\Phi$ and $\nu_i$ as shown in Figure~\ref{slopeconvention}.  By Proposition~\ref{translation1} the slope of $\Phi$ is $\frac{1}{c_i-n_i}$ while the slope of $\nu_i$ is $-\frac{1}{n_i}$.  The form of the intervals given, depends on a case by case analysis of whether or not they contain slope $\infty$.  Weak symplectic fillability   and universal tightness   of $(\xi_+)_{(S,h)}$ follows from Proposition~\ref{oppfoln}.
\end{proof}

One can translate the results of \cite{HKM2,R2} into the current context as follows.

\begin{thm}\label{existence} \cite{R2} When the binding $L$ is connected, $c>0$, and the monodromy $h$ has pseudo-Anosov representative, there are $\partial M_0$-taut foliations in $M_0$ transverse to $\Phi_0$ which realize all slopes in an interval  $J$ as follows:
\begin{enumerate}
\item $c=n$ and $J=(-\infty,\infty)$, or 
\item $c>n$ and $J=(-\infty,1)$, or
\item $c<n$ and $J=(-1,\infty)$.
\end{enumerate}
\end{thm}

The next corollary follows by intersecting the intervals where foliations exist in Theorem~\ref{existence} with the intervals where they are required in Proposition~\ref{oppfolnTh}.

\begin{cor} There exists $\mathcal F_-$ as described in Proposition~\ref{oppfolnTh}, if one of the following is true
\begin{enumerate}
\item $c=n$ and $(-\infty,-\frac{1}{n})\cap (-\infty,\infty)\ne\emptyset$, or
\item $c>n$ and $((\frac{1}{c-n},\infty] \cup [-\infty,-\frac{1}{n}))\cap (-\infty,1)\ne\emptyset$, or
\item $c<n$ and $ (-\frac{1}{n-c},-\frac{1}{n})\cap (-1,\infty)\ne\emptyset$. \qed
\end{enumerate}
\end{cor}

The case when the fractional Dehn twist coefficient is greater than or equal to $1$ is of particular interest since

\begin{cor} If $c \ge 1$ there exists $\mathcal F_-$ as described in Proposition~\ref{oppfolnTh}.
\end{cor}

\begin{proof}
If $c>0$ then $n\ge 0$, with $n=0$ only when $c\in (0,1/2]$.  Thus $c \ge 1$ implies $n\ge 1$, and it follows that the intersections in Cases (1) and (2) are nonempty.  In Case (3), $1 \le c < n$, and the intersection is again nonempty.
\end{proof}

This corollary is exactly what is needed to complete the proof of a theorem of Honda, Kazez, and Mati\'c in \cite{HKM2} that was one of the original motivations for this work.

\begin{thm} If $(S,h)$ is an open book decomposition such that $S$ has connected boundary, $h$ is isotopic to a pseudo-Anosov homeomorphism, and the fractional Dehn twist coefficient of $h$ is greater than or equal to $1$, then the contact structure canonically associated to the open book decomposition, $\xi(S,h)$, is weakly symplectically fillable   and universally tight.
\end{thm}

The proof strategy of \cite{HKM2} used a single foliation $\mathcal F$  defined on all of $(S,h)$ as constructed by Roberts, \cite{R1,R2} with boundary slope related to open book data.  Next they wanted to apply the Eliashberg-Thurston theorem to produce a weakly symplectically fillable  and universally tight contact structure $\xi$.  Finally they argued that the two contact structures $\xi(S,h)$ and $\xi$ were necessarily equivalent.

There is a somewhat surprising aspect that arises in addressing the issue of lack of smoothness of $\mathcal F$. It is that we do not use the same foliation.  To complete their proof with our strategy, we require the existence of two foliations $\mathcal F^+$ and $\mathcal F^-$, both of which exist   on the complement of the binding of $(S,h)$, by the work of Roberts, and have boundary slopes on either side of the boundary slope of $\mathcal F$.  One foliation is used to the produce a positive contact structure, the other a negative contact structure, and both are necessary to conclude that the approximating contact structure is weakly symplectically fillable   and universally tight.

The notion of $\Phi$-approximating contact structure that we produce is sufficient to conclude that $\xi(S,h)$ and $\xi$ are equivalent using the argument of \cite{HKM2}. \qed
\medskip

When $M=(S,h)$ has binding which is not connected, our results can be applied to the Kalelkar-Roberts constructions  of $\partial M_0$-taut foliations in $M_0$ transverse to $\Phi_0$.

\begin{thm} [Theorem~1.1,\cite{KRo}] There are $\partial M_0$-taut oriented foliations in $M_0$ transverse to $\Phi_0$ and realizing a neighborhood of rational boundary multislopes about the boundary multislope of the fibration.
\end{thm}

\begin{cor} Suppose $M$ has open book decomposition $(S,h)$ and fractional Dehn twist coefficients $c_i, 1\le i\le b$. There are constants $A_i= A_i(M_0), 1\le i\le b$, dependent on $M_0$ such that if $c_i>A_i$, then $\xi_{(S,h)}$ is weakly symplectically fillable   and universally tight.
\end{cor}

Work of  Baldwin and Etnyre \cite{BE} implies that any such constants $A_i$ must depend on $M_0$, at least in the case that the page $S$ has genus one.

\begin{thm}[Theorem 1.16 of \cite{BE}]
There exist open books whose fractional Dehn twist coefficients are arbitrarily large, but whose compatible contact structures are not $C^0$ close to  smooth orientable taut foliations.
\end{thm}

\begin{cor}
There exist open books whose fractional Dehn twist coefficients are arbitrarily large, but whose compatible contact structures are not $\Phi$-close to taut oriented bracketed $C^{1,0}$-foliations, for any volume preserving flow $\Phi$.
\end{cor}

\appendix
\section{Some symplectic topology}\label{symplectic}

This section contains an overview of the relationship between foliations, volume preserving flows, symplectic topology, and contact topology that is summarized in Theorem~\ref{weaklysymplectic}.

Let $\mathcal F$ be a transversely oriented, taut $C^{1,0}$-foliation in $M$.  Fix a metric on $M$, and let $\Phi$ be a volume preserving flow transverse to $\mathcal F$.  The starting point for the interconnections we will describe is a carefully chosen 2-form.

\begin{definition} (See Section~3.2 of \cite{ET}.)
Let $\xi$ be a co-oriented $C^k$ 2-plane field on a smooth 3-manifold $M$ with $k\ge 0$. A smooth closed 2-form $\omega$ on $M$ is said to {\it dominate} $\xi$ if $\omega|\xi$ does not vanish (i.e., if  $p\in M$ and $X_p,Y_p$ is a basis for $\xi_p$, then $\omega_p(X_p,Y_p)\ne 0$). A smooth closed 2-form $\omega$ on $M$ is said to {\sl positively dominate} $\xi$ if $\omega|\xi$ is positive (i.e., for all $p\in M$, if $X_p,Y_p$ is a positively oriented basis for $\xi_p$, then $\omega_p(X_p,Y_p)> 0$).
\end{definition}

To produce such a dominating closed 2-form, let $\Omega$ be the volume form on $M$ preserved by the smooth flow $\Phi$, and let $X$ be the vector field which generates $\Phi$.  Define
$$\omega = X \lrcorner\Omega.$$
Recall that $\Phi$ is volume preserving if and only if $\mathcal L_X\Omega = 0$, where $\mathcal L_X$ denotes the Lie derivative with respect to $X$. (See for example, Proposition 18.16 of \cite{Lee}.) By Cartan's Formula (see, for example, Proposition 18.13 of \cite{Lee}),
$$\mathcal L_X\Omega = X\lrcorner (d\Omega) + d(X\lrcorner \Omega) = d(X\lrcorner \Omega) = d\omega.$$
It follows that $\omega$ is closed.  From its definition, $\omega$ is killed by the flow direction $X$, thus a co-oriented 2-plane field $\xi$ is positively dominated by $\omega$ if and only if it is everywhere positively transverse to $\Phi$.

A closed 2-form dominating $T\mathcal F$ can be produced directly from a taut foliation \cite{Sullivan,Hass}  thereby eliminating the need to choose a metric, a volume form, and a flow preserving the volume form.  We have chosen to emphasize the volume preserving flow since it clarifies the local nature of our foliations in flow boxes.

Specialize now to the case of Theorem~\ref{weaklysymplectic} in which $\xi$ is a contact structure positively transverse to $\Phi$.  Choose a 1-form $\alpha$ such that ker$\,\alpha=\xi$ and $\alpha \wedge \omega>0$.  Define a 2-form $\tilde{\omega}$ on $M \times [-1,1]$ using the projection map $p$ and the formula

\begin{equation}\label{sympform}
\tilde{\omega} = p^{\star}(\omega) + \epsilon d(t\alpha).
\end{equation}

Direct computation shows that if $\epsilon$ is positive and small enough, $(M\times[-1,1],\tilde{\omega})$ is a symplectic manifold with boundary, that is, $\tilde{\omega}$ is a non-degenerate 2-form.  The important role of the positive and negative contact structures, $\xi_+$ and $\xi_-$, of Theorem~\ref{weaklysymplectic} will be described after the next definition.

\begin{definition} A  boundary component $Y$ of a symplectic manifold $(W,\tilde{\omega})$ is called {\sl weakly convex} if $Y$ admits a positive contact structure dominated by $\tilde{\omega}|_Y$.
\end{definition}

This is precisely the structure that the boundary components of  $(M\times[-1,1],\tilde{\omega})$ have. The restriction of $\tilde{\omega}$ positively dominates $\xi_+$ on $M\times\{1\}$.  Because of boundary orientations, $\xi_-$ defines a positive contact structure on $M\times\{-1\}$ that is positively dominated by the restriction of $\tilde{\omega}$.  Moreover, $\tilde{\omega}$ restricts to $\omega$ on $M\times\{0\}$, thus $M \times[-1,0]$ also has weakly convex boundary.

Since  both  boundary components of either $(M\times[-1,1], \tilde{\omega})$ or the restriction of $\tilde{\omega}$ to $M\times [-1,0]$ are weakly convex, they give examples of weak symplectic fillings.

\begin{definition} A {\sl weak symplectic filling} of a contact manifold $(M,\xi)$ is a symplectic manifold $(W,\tilde{\omega})$ with $\partial W = M$ (as oriented manifolds) such that $\tilde{\omega}|_{\xi}> 0$. A contact manifold $(M,\xi)$ which admits a weak symplectic filling is called {\sl weakly symplectically fillable}.
\end{definition}

\begin{definition} 
A {\sl strong symplectic filling} of a contact manifold $(M,\xi)$ is a symplectic manifold $(W,\tilde{\omega})$ with $\partial W = M$ (as oriented manifolds) where $\xi = \text{ker} \, \alpha$ for a 1-form $\alpha$ satisfying $d\alpha = \tilde{\omega}|_M$. A contact manifold $(M,\xi)$ which admits a strong symplectic filling is called {\sl strongly symplectically fillable}.
\end{definition} 

In general, strong symplectic fillability is a stronger condition than weak symplectic fillability \cite{Eli2}. However, by Lemma 1.1 of \cite{OO}, when $M$ is a rational homology sphere, $(M,\xi)$ is weakly symplectically fillable if and only if it is strongly symplectically fillable.

\begin{definition}
A contact manifold is said to be weakly (or strongly) {\sl semi-fillable} if it is one boundary  component of a weak (or strong)  symplectic filling.
\end{definition}

A weakly (or strongly) semi-fillable contact manifold is weakly (or strongly) fillable \cite{Eli, Et}.

The following fundamental theorem gives an example of the importance of weak symplectic fillability in contact topology.

\begin{thm} [\cite{Gr},\cite{Eli1},\cite{ET}] \label{EGr} Weakly symplectically fillable contact structures are  tight.   Weakly symplectically fillable contact structures that are $\Phi$-close to taut foliations are universally tight.
\end{thm}

\end{document}